\documentclass[journal,10pt]{IEEEtran}
\newcommand{\secref}[1]{Section~\ref{#1}}
\newcommand{\figref}[1]{Fig.~\ref{#1}}
\newcommand{\tabref}[1]{Table~\ref{#1}}

\usepackage{pifont}
\usepackage{makecell}
\usepackage{markdown}
\usepackage{tcolorbox}
\usepackage{amsmath}
\usepackage{graphicx}
\usepackage{subcaption}
\usepackage{makecell}

\usepackage{graphicx}
\usepackage{amsmath}
\usepackage{comment}
\usepackage{bm}
\usepackage{multirow}
\usepackage{threeparttable,booktabs}
\usepackage{tikz}
\usepackage{balance}
\usepackage{courier}                                       
\usepackage[mathcal]{eucal}
\usepackage[noend]{algpseudocode}
\algrenewcommand\textproc{\texttt}
\makeatletter\let\float@addtolists\relax\makeatother
\usepackage{algorithm}
\usepackage{filecontents}                                  
\usepackage{pgfplots}
\usepackage{pgfplotstable}
\usepackage{algpseudocode}                                 
\usepackage{algorithm}
\usepackage{graphicx}                                      
\usepackage{amsmath}
\usepackage{amssymb}
\usepackage{amsfonts}
\usepackage{booktabs}
\usepackage{multirow}
\usepackage[square, comma, sort&compress, numbers]{natbib}
\usepackage{comment}
\usepackage{soul}                                          
\soulregister\cite7
\soulregister\ref7
\soulregister\pageref7
\usepackage{amsthm}
\usepackage{etoolbox}                                      
\usepackage{url}
\usepackage{nth}                                           
\usepackage{bm}                                            
\usepackage{courier}
\usepackage{balance}
\usepackage{threeparttable}
\usepackage[bookmarks=false]{hyperref}                     %
\hypersetup{
    colorlinks = true,
    citecolor  = blue,
    linkcolor  = blue,
    urlcolor   = blue,
}
\usepackage{filecontents}                                  
\usepackage{tikz}
\usetikzlibrary{positioning, fit}
\usepackage{pgfplots}
\usepackage{pgfplotstable}
\pgfplotsset{compat=newest}
\usepackage{caption} 
    
\usepackage{tabu}
\usepackage{tabularx}
\usepackage{rotating}
\usepackage{supertabular}
\usepackage{tikz}
\usepackage{hhline}
\usepackage{longtable}
\usepackage{svg}
\usepackage{rotating}
\usepackage{xcolor}
\usepackage{bbding}

\usepackage{color}
\usepackage{xcolor}
\usepackage{threeparttable}
\usepackage{utfsym}
\usepackage{url}
\usepackage{hyperref}
\usepackage{cleveref} 
\usepackage{multirow}
\usepackage{amsmath}
\usepackage{float}
\usepackage{algorithm}
\usepackage[export]{adjustbox}
\usepackage[tableposition = top]{caption} 
\usepackage{booktabs}                     
\usepackage[per-mode = symbol]{siunitx}   
\usepackage{colortbl,booktabs}
\usepackage{multicol}
\usepackage{multirow} 
\usepackage{multirow, makecell, caption}
\usepackage{diagbox}

\usepackage{amsmath,amssymb,amsfonts}
\usepackage{graphicx}
\usepackage{textcomp}
\usepackage{xcolor}

\usepackage{lettrine}

\graphicspath{{./Images/}}

\definecolor{CUHKorange}{RGB}{244,106,18} 
\definecolor{CUHKblue}{RGB}{0,111,190}    
\definecolor{CUHKgreen}{RGB}{0,127,128}   
\definecolor{CUHKred}{RGB}{228,46,36}     
\definecolor{CUHKyellow}{RGB}{198,148,34} 
\definecolor{CUHKdark}{RGB}{114,44,114}   
\definecolor{CUHKmiddle}{RGB}{144,44,144} 
\definecolor{CUHKlight}{RGB}{167,44,167} 

\algrenewcommand\textproc{\texttt}

\makeatletter
\let\OldStatex\Statex
\renewcommand{\Statex}[1][3]{%
  \setlength\@tempdima{\algorithmicindent}%
  \OldStatex\hskip\dimexpr#1\@tempdima\relax
}
\makeatother

\RequirePackage[normalem]{ulem} 
\RequirePackage{color}\definecolor{RED}{rgb}{1,0,0}\definecolor{BLUE}{rgb}{0,0,1} 

\graphicspath{{./figs/}}

\title{
\texttt{ATMPlace}: Analytical Thermo-Mechanical-Aware Placement Framework for 2.5D-IC
}

\author{
    Qipan Wang~\IEEEmembership{Student Member,~IEEE}, Tianxiang Zhu~\IEEEmembership{Student Member,~IEEE}, \\
    Tianyu Jia~\IEEEmembership{Member,~IEEE}, 
    Yibo Lin~\IEEEmembership{Member,~IEEE}, 
    Runsheng Wang~\IEEEmembership{Member,~IEEE}, 
    Ru Huang~\IEEEmembership{Fellow,~IEEE}
    \thanks{ 
        This work was supported in part by the National Science Foundations of China (Grant No. 62125401, 62034007), the Natural Science Foundation of Beijing, China (Grant No. Z230002), the Grant QYJS-2023-2303-B, the Beijing Outstanding Young Scientist Program (JWZQ20240101004), and the 111 project (B18001).
        The preliminary version has been presented at the 43rd IEEE/ACM International Conference on Computer-Aided Design in 2024~\cite{wang2024atplace2}.
    }
    \thanks{Q.~Wang is with the School of Advanced Interdisciplinary and the School of Integrated Circuits, Peking University, Beijing, China.}
    \thanks{T.~Zhu and T.~Jia are with the School of Integrated Circuits, Peking University, Beijing, China.}
    \thanks{Y.~Lin, R.~Wang, and R.~Huang are with the School of Integrated Circuits, Peking University, Beijing Advanced Innovation Center for Integrated Circuits, Beijing, China. and Institute of Electronic Design Automation, Peking University, Wuxi, China}
    \thanks{Corresponding authors: Yibo Lin (\url{yibolin@pku.edu.cn}), and Runsheng Wang (\url{r.wang@pku.edu.cn}}
    
}

\begin{document}
\maketitle
\thispagestyle{empty}
\begin{abstract}
Rising demand in AI and automotive applications is accelerating 2.5D-IC adoption, with multiple chiplets tightly placed to enable high-speed interconnects and heterogeneous integration. As chiplet counts grow, traditional placement tools—limited by poor scalability and reliance on slow simulations—must evolve beyond wirelength minimization to address thermal and mechanical reliability, critical challenges in heterogeneous integration.

In this paper, we present \texttt{ATMPlace}, the first analytical placer for 2.5D-ICs that jointly optimizes wirelength, peak temperature, and operational warpage using physics-based compact models. 
It generates Pareto-optimal placements for systems with dozens of chiplets. Experimental results demonstrate superior performance: 146\% and 52\% geo-mean wirelength improvement over \texttt{TAP-2.5D} and \texttt{TACPlace}, respectively, with 3–13\% lower temperature and 5–27\% less warpage — all achieved $\sim$10$\times$ faster. 
The proposed framework is general and can be extended to enable fast, scalable, and reliable design exploration for next-generation 2.5D systems.
\end{abstract}

\begin{IEEEkeywords}
2.5D-IC; Thermo-Mechanical Optimization; Chiplet Placement; Thermal Warpage
\end{IEEEkeywords}

\section{Introduction} \label{sec:Introduction}
\begin{figure}[tbh]
    \includegraphics[width=1.0\linewidth]{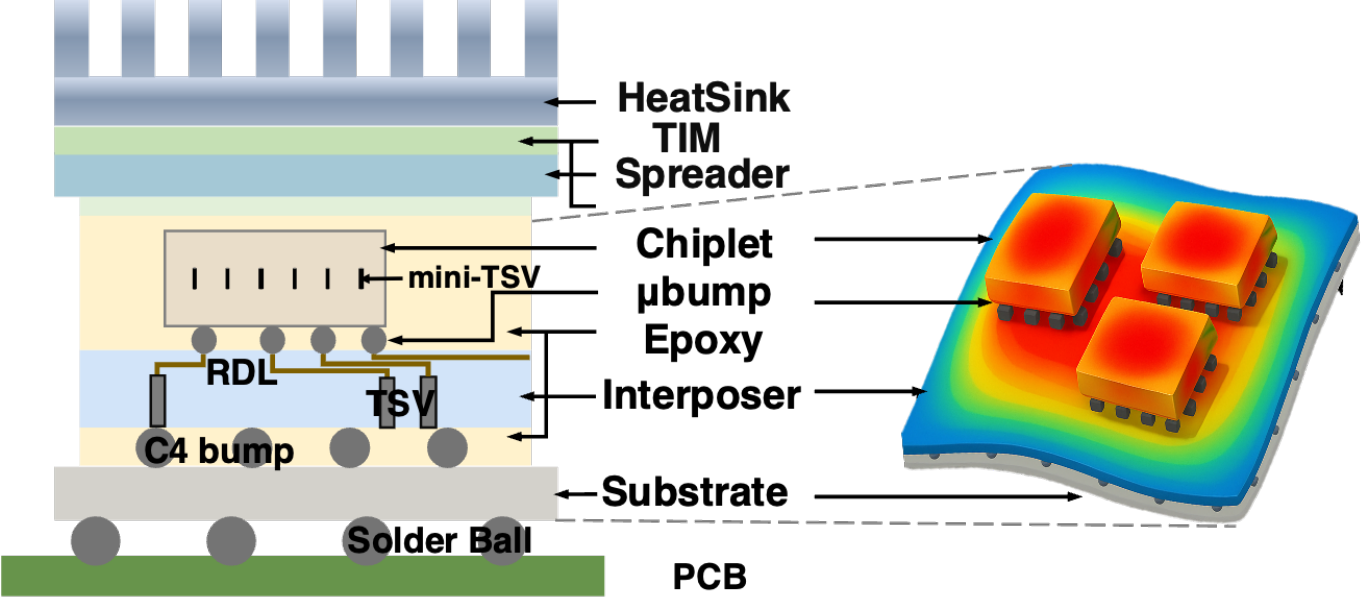}
    \caption{Layered cross-section (left) and isometric illustration of thermal deformation (right) in 2.5D-IC structure.} 
    \label{fig:25dic}
\end{figure}

\begin{table}[tbh]
\centering
\caption{Comparison of features of different chiplet placement studies. \textbf{Scalability} represents the maximum number of chiplets the method can handle.}
\footnotesize
\resizebox{0.49\textwidth}{!}{
\begin{tabular}{|c|c|cc|cc|}
\hline
    \textbf{Algorithm} & \textbf{Work} & \textbf{Wirelength} & \textbf{Scalability} & \thead{Thermal\\aware} & \thead{Warpage\\aware} \\ \hline
        \multirow{4}{*}{Heuristic} 
        & \citet{ho2013multiple} & \multirow{4}{*}{Medium} & Low & \ding{55} & \ding{55} \\
        & \texttt{TAP-2.5D}~\cite{ma2021tap} &  & Low & \ding{51} & \ding{55} \\ 
        & \texttt{TACPlace}~\cite{yu2025tacplace} & & Medium & \ding{51} & \ding{55} \\ 
        & \citet{10680720} & & Medium & \ding{51} & \ding{55} \\ \hline
        \multirow{2}{*}{Enumeration} 
        & \citet{osmolovskyi2018optimal} & \multirow{2}{*}{\textbf{Low}} & Medium & \ding{55} & \ding{55} \\
        & \texttt{SP-CP}~\cite{chiou2023chiplet} & & Medium & \ding{51} & \ding{55} \\ \hline
        \multirow{2}{*}{RL} 
        & \citet{duan2023rlplanner} & \multirow{2}{*}{Medium} & Medium & \ding{51} & \ding{55} \\
        & \citet{deng2024chiplet} &  & Medium & \ding{51} & \ding{55} \\ \hline
        \multirow{3}{*}{Analytical} 
        & \citet{chen2023floorplet} & \multirow{3}{*}{\textbf{Low}} & Medium & \ding{55} & \ding{51} \\
        & \texttt{ATPlace2.5D}~\cite{wang2024atplace2} & & \textbf{High} & \ding{51} & \ding{55} \\
        & \texttt{ATMPlace} (ours) & & \textbf{High} & \ding{51} & \ding{51} \\ \hline
    \end{tabular}}
    \label{table:comp}
\end{table}

The continued slowdown of technology scaling has rendered monolithic SoCs increasingly cost-prohibitive, especially in domains that demands both high performance and rapid scalability, such as automotive electronics and AI accelerators. 
To circumvent these limitations, 2.5D integration has emerged as a practical path forward, assembling pre-verified chiplets on an interposer that delivers high-bandwidth communication, technology heterogeneity, and cost reduction~\cite{naffziger2021pioneering}. By enabling reuse and modularity, 2.5D-IC design promises to reshape the semiconductor design paradigm.

However, the effectiveness of 2.5D-ICs fundamentally depends on the placement (floorplan) problem~\cite{Survey_of_25D}.\footnote{We use the terminology “placement’’ and “floorplan’’ interchangeably to denote the physical arrangement of chiplets.}
In this stage, the positions and orientations of heterogeneous chiplets must be determined under multiple performance, thermal, and manufacturability constraints.

A wide spectrum of methods has been explored, including heuristic methods (simulated annealing (SA)~\cite{ho2013multiple,coskun2018cross,ma2021tap,parekh2025stamp25dstructuralthermalaware}, genetic algorithm (GA)~\cite{yu2025tacplace,10680720}), enumeration-based search~\cite{liu2014floorplanning,chiou2023chiplet}, reinforcement learning (RL)~\cite{duan2023rlplanner,deng2024chiplet}, and analytical optimization~\cite{wang2024atplace2,romano2025diffchip}. 
Heuristic methods rely on multi-objective cost functions to co-optimize metrics. However, they require long runtime and often converge to suboptimal layouts.
Enumeration-based strategies can deliver good solutions for designs of small instances. Yet they become infeasible when the number of chiplets exceeds ten due to exponential complexity.
RL-based approaches offer a data-driven alternative, but they need extensive training data, causing unstable performance.
All three categories exhibit limited scalability as chiplet count and diversity increase. For example, heuristic and enumeration methods may require hours for designs with $\geq 10$ chiplets, which is impractical for iterative design-space exploration~\cite{chen2023floorplet}.
In contrast, analytical approaches can handle configurations with up to 60 chiplets in several minutes, as demonstrated in~\cite{wang2024atplace2}.

Moreover, prior work primarily tries to minimize area and wirelength, frequently yielding compact placements with high power density and thermal failure risks. 
Recent thermal-aware placement studies either embed temperature constraints ~\cite{coskun2018cross}, augment cost functions with worst-case temperature~\cite{ma2021tap}, or perform post-optimization refinement~\cite{chiou2023chiplet}. 
These approaches rely on iterative numerical thermal simulations~\cite{stan2003hotspot,6241575}, significantly inflating runtime. Treating temperature as a hard constraint (e.g., limiting peak temperature or chiplet proximity~\cite{ma2021tap,chiou2023chiplet}) also narrows the exploration of the broader thermo-design space required by large systems.

Finally, warpage poses a critical reliability challenge in large 2.5D packages with a thin interposer. Industrial reports have highlighted warpage issues in advanced GPU productions~\cite{SemiAnalysis2024Blackwell}. 
While prior studies analyze and optimize fabrication-stage warpage~\cite{hsu2022transitive,chen2023floorplet}, operational warpage driven by thermal-expansion mismatch and thermal gradients during operation is largely overlooked. 
The analytical model in~\cite{9503285} targets 2D-ICs but does not capture the geometric and material complexity of 2.5D assemblies.

To conclude, despite a decade of research, current 2.5D-IC placement algorithms remain far from meeting industrial needs. We summarize the features and abilities of previous works in \tabref{table:comp}. In this work, we propose an analytical placement framework for 2.5D-ICs with joint thermo-mechanical optimization. To overcome scalability limits, we adopt an orientation-aware analytical placer and couple it with fast physics-based compact models of temperature and warpage. Our main contributions are:
\begin{itemize}
  \item \texttt{ATMPlace}: an analytical chiplet placement framework that jointly optimizes total wirelength and \emph{thermo-mechanical reliability} (temperature and warpage) within a unified flow.
  \item Two mixed-integer linear programming (MILP) formulations for initialization and legalization that assign chiplet orientations and reduce wirelength while preserving feasibility.
  \item A physics-based analytical compact model for fast thermal and mechanical evaluation integrated into the optimization loop. It achieves a mean correlation $>0.97$ and a $>8000\times$ speedup for both thermal and mechanical evaluation versus the numerical simulator.
\end{itemize}

We evaluate our method on large-scale 2.5D-ICs. Under thermo-mechanical-aware placement, it outperforms \texttt{TAP-2.5D} and \texttt{TACPlace} simultaneously in wirelength (geometric mean improvement: 146\% and 52\%), peak temperature (reduction: 3\% and 13\%), and warpage (reduction: 5\% and 27\%), with $\sim$10$\times$ speedup. In wirelength-driven mode, it matches \texttt{SP-CP}'s total wirelength within 1\% deviation.  
The rest of the paper is organized as follows: \secref{sec:Preliminary} introduces background; \secref{sec:Algorithm} details the algorithm; \secref{sec:Results} presents results; \secref{sec:Conclusion} concludes the paper.

\begin{figure}[tbh]
    \includegraphics[width=0.95\linewidth]{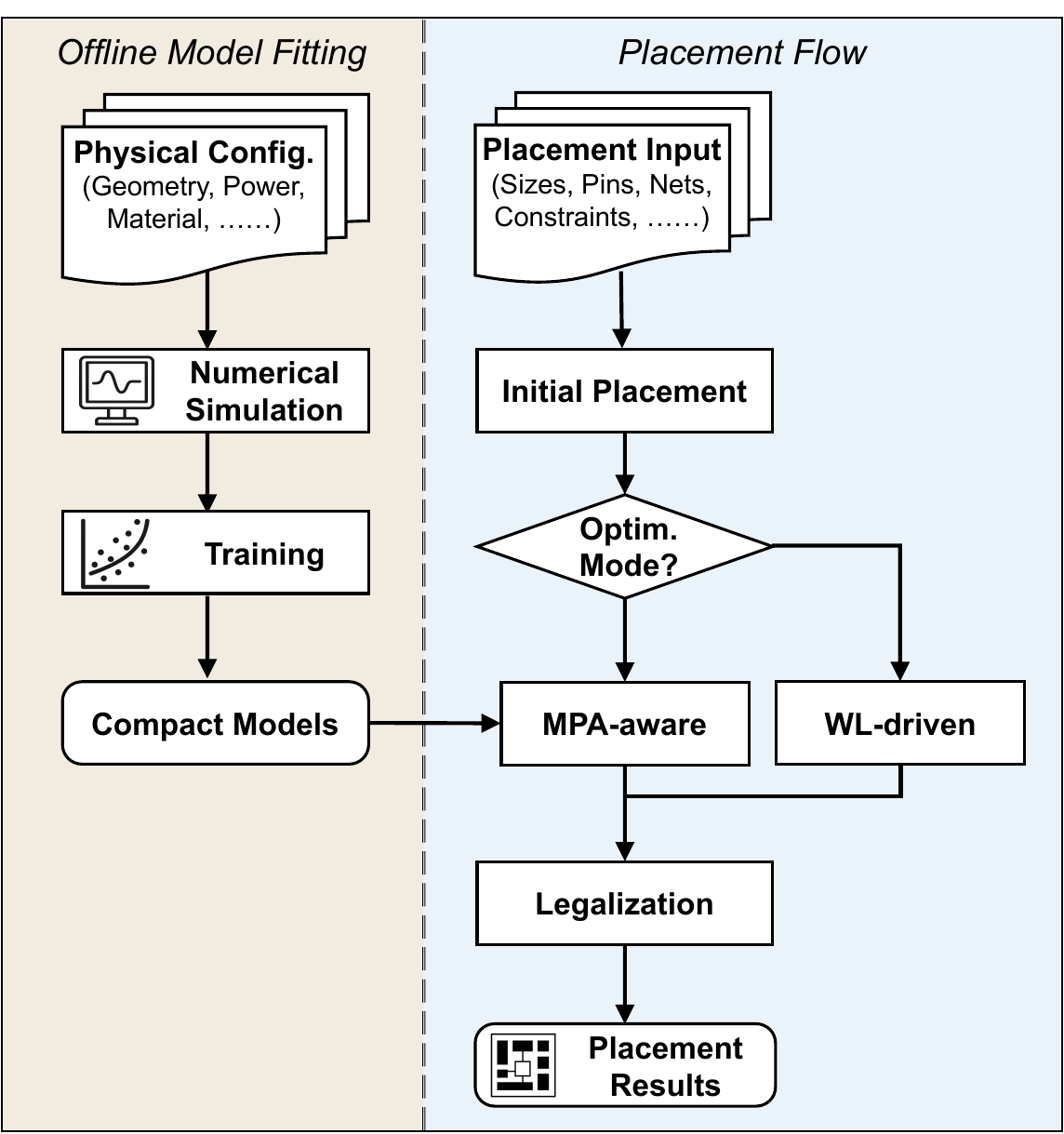}
    \caption{Illustration of the \texttt{ATMPlace} framework, including the compact model fitting and placement flow.}
    \label{fig:framework}
\end{figure}

\begin{figure}[tbh]
    \centering
    \includegraphics[width=0.95\linewidth]{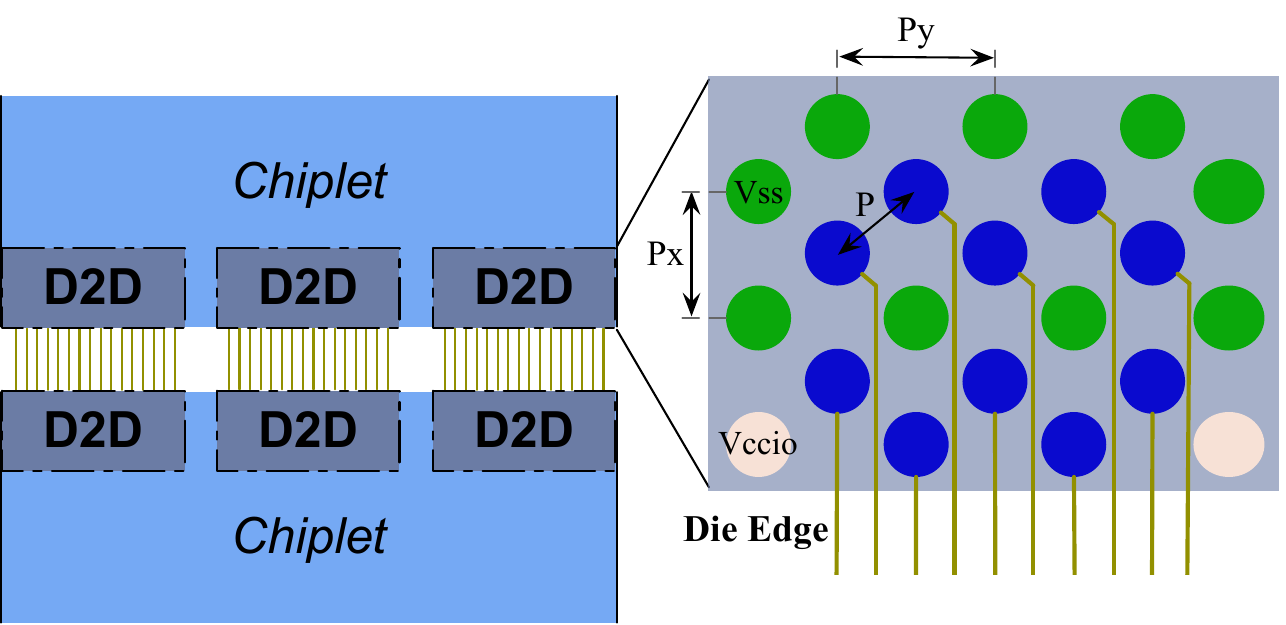}
    \caption{Illustration of the D2D link between chiplets (left), and a portion of the detailed package bump map of a $\times16$ module (right). The seashell and green circles are the I/O supply (Vccio) and ground reference (Vss) bump, while the blue circles represent the $Tx/Rx$ bumps, respectively.}
    \label{fig:bump}
\end{figure}

\section{Preliminaries}
\label{sec:Preliminary}
In this section, we begin by describing the 2.5D-IC architecture considered in this work in~\secref{sec:25dconfig}. Then we present the thermal and warpage models in~\secref{sec:thermal} and~\secref{sec:warpage-sim}. Finally, \secref{sec:metric} formally defines the placement optimization objectives and constraints.

\subsection{2.5D-IC Configuration} \label{sec:25dconfig}

The simplified 2.5D-IC structure studied in this work is illustrated in~\figref{fig:25dic}. It consists of a passive interposer mounted on an organic substrate via C4 bumps. 
We adopt a transistor-free passive interposer due to its lower cost and design flexibility, while the proposed framework can be extended to active variants. 
Through-silicon vias (TSVs) pass through the interposer to connect chiplets to external circuits.
Chiplets sit atop the interposer and are modeled as silicon blocks with given power dissipation. They interface with the redistribution layers (RDL) through microbump arrays. Due to their small size, these microbumps are treated as a homogeneous layer with an effective thermal conductivity. 
The regions between chiplets and microbumps are filled with epoxy, introducing material heterogeneity across the package.

We assume a library of pre-fabricated, modular chiplets is available. Each chiplet includes one or more interface units for data exchange. These interfaces connect to the interposer via microbumps, enabling communication between chiplets and forming complete systems for various applications.
To support interoperability among diverse chiplets, standardized interfaces are essential for building a robust chiplet ecosystem. Two primary types of die-to-die (D2D) interfaces exist: serial and parallel. Serial links use a few differential pairs for data transmission, while parallel interfaces employ tens to hundreds of connections. These differ in bandwidth, latency, and packaging requirements, making universal integration across heterogeneous technologies a complex challenge. 

For demonstration purposes, we adopt the Universal Chiplet Interconnect Express (UCIe) standard~\cite{sharma2022universal} for D2D communication. This choice enables meaningful evaluation and supports future research, though our placement framework is compatible with other interface standards. As shown in the left side of \figref{fig:bump}, D2D links connect chiplets through microbump pairs. The right side shows a UCIe reference bump matrix: each signal lane consists of a transmitter (Tx) and receiver (Rx) pair, connected across chiplets. For complex topologies, a hub chiplet may be used to route traffic among multiple chiplets~\cite{yang2023challenges}.

\subsection{Thermal Model} \label{sec:thermal}
We analyze steady-state temperature under the worst-case power profile—a standard choice for thermal-aware design.  
The governing heat conduction equation is:
\begin{equation}
    \nabla\cdot\left(\kappa(\textbf{r})\nabla T(\textbf{r})\right) = -\textbf{P}(\textbf{r}), \label{equ:thermal}
\end{equation}
subject to appropriate boundary conditions.
Here, $\kappa(\mathbf{r})$ denotes the spatially varying (and potentially anisotropic) thermal conductivity, and $\mathbf{P}(\mathbf{r})$ is the volumetric power density distribution of the chiplets. Our thermal model naturally supports multiple workload-dependent power profiles: handling $M$ scenarios only requires $M$ independent evaluations to track the worst-case response, without modifying the solver.  

Following the common practice~\cite{ma2021tap,chiou2023chiplet}, we model heat dissipation through a primary path: chiplet $\rightarrow$ TIM $\rightarrow$ heat spreader $\rightarrow$ heatsink $\rightarrow$ ambient air. A secondary path via substrate and PCB is also modeled, and lateral surfaces are assumed adiabatic. We use an air-cooled heatsink with convective resistance of $0.1\mathrm{K/W}$~\cite{6241575}. 
Layer sizes, material properties, bump/TSV dimensions follow previous works~\cite{ma2021tap,wang2024atsim}.

Various thermal simulators have been proposed, and mainstream methods encompass numerical and analytical approaches \cite{pedram2006thermal}. Numerical methods that mesh the system and solve linear equations form the foundation of a majority of simulators.
However, they also demand substantial computation time for high precision, and are hard to derive the gradient information, apart from the adjoint method which brings about high runtime cost. 
Hence, analytical methods are favored during the design stage, where iterative simulations are required. 
Such methods can provide rapid solutions through closed-form approximate expressions, either obtaining exact solutions based on simplified models~\cite{wang2007accelerated} or constructing approximate expressions based on accurate simulation results~\cite{Torregiani2009compact,kung2011thermal,ziabari2014power}. The latter approach often relies on Green's function, which is based on the impulse response of the system concerning input power given specified initial or boundary conditions.

\subsection{Mechanical Model} \label{sec:warpage-sim}
Thermal gradients induce mechanical warpage due to the coefficient of thermal expansion (CTE) mismatch across heterogeneous materials (e.g., Si dies, organic substrates, Cu pillars). 
To assess thermally induced mechanical deformation in the heterogeneous integration stack, we model the out-of-plane displacement $w(\mathbf{r})$ governed by the biharmonic plate equation under non-uniform thermal loading~\cite{lo2024efficient}:
\begin{equation}
    \nabla^4 w(\mathbf{r}) = \frac{1 - \nu}{E h^3} \nabla^2 \big( \alpha h \, \Delta T(\mathbf{r}) \big), \label{equ:warpage}
\end{equation}
where $E$, $\nu$, and $\alpha$ denote the effective Young’s modulus, Poisson’s ratio, and coefficient of thermal expansion of the equivalent laminate stack, respectively; $h$ is the effective thickness, and $\Delta T(\mathbf{r}) = T(\mathbf{r}) - T_{\text{ref}}$ represents the local temperature rise obtained from the thermal model in~\secref{sec:thermal}. Boundary conditions are assumed to be rigid motion elimination at the package edges as previous~\cite{Suhir1989,lo2024efficient}, and interfacial stresses are captured through homogenized material properties.

The warpage metric is critical for assembly yield and reliability, as it quantifies the global mechanical deformation of the package. It is defined as the peak-to-valley deviation of the out-of-plane displacement field:
    $\max_{\mathbf{r}} w(\mathbf{r}) - \min_{\mathbf{r}} w(\mathbf{r})$.
This formulation enables efficient co-analysis of runtime thermal and mechanical effects.

To simplify finite element analysis, the analytical model proposed in \cite{Suhir1989} has been widely adopted for its simplicity and computational efficiency on simple bi-layer or tri-layer structures. However, its accuracy is limited when applied to modern heterogeneous packaging with complex material stacks and non-uniform thermal profiles. More recently, \citet{tasi2020} proposed a warpage model incorporating polynomial expressions with assumed curvature, which is more suitable for warpage estimation in advanced packaging. 

\begin{table}[htbp]
    \centering
    \caption{Some notations used in this paper.}
    \begin{tabular}{|c|p{0.75\linewidth}|}
        \hline 
        \multicolumn{1}{|c|}{\textbf{Notation}} & \multicolumn{1}{c|}{\textbf{Meaning}}  \\
        \hline 
        $\mathbb{C, A, E}$ & Sets of chiplets, connected chiplet pairs, and nets. \\\hline
        $i,j; C_i,C_j$ & Index and symbol of a chiplet $\in \mathbb{C}$. \\
        \hline
        $x_i, y_i, \theta_i$ & X and Y coordinates of the center, and rotation angle of the chiplet $C_i$. \\
        \hline
        $w_i, h_i, t_i$ & Width, height, and thickness of chiplet $C_i$. \\
        \hline
        $w_i^\prime, h_i^\prime$ & Width and height of chiplet $C_i$ after rotation. \\
        \hline
        $\mathbb{B}_i, \mathbb{B}_e$ & Set of microbumps belongs to the chiplet $C_i$ or net $e$. \\
        \hline
        $x_{p_i}, y_{p_i}$ & X- and Y-offsets from the center of the chiplet $C_i$ for a microbump $p_i\in\mathbb{B}_i$. \\
        \hline
        $w_{\text{gap}}$ & Minimum spacing between two chiplets, set to be $100 \mu m$ as in \cite{ma2021tap}. \\
        \hline
        $W, H$ & Width and height of the placement region (interposer). \\
        \hline
    \end{tabular}
    \label{tab:notat}
\end{table}

\subsection{Problem Formulation} \label{sec:metric}
We list key notations in \tabref{tab:notat}. The goal is to determine the position and orientation (counterclockwise direction) of the chiplets, minimizing the total wirelength, with runtime temperature and warpage satisfying certain conditions meanwhile. 
Note that the orientation must adhere to certain legal values, denoted by $\Theta=\{\theta^0,\theta^1,\theta^2,\theta^3\},$ where $\theta^{0\sim3}$ are ${0^\circ,90^\circ,180^\circ,270^\circ}$ respectively. 

In this work, inter-chiplet communication cost is quantified by total wirelength, consistent with prior work~\cite{liu2014floorplanning,ma2021tap}. 
We focus on signal nets (\textit{i.e.}, Tx–Rx pairs in D2D links), assuming power, ground, and I/O connections are handled later. Since these nets are typically two-pin and routed as transmission lines in the interposer, we avoid traditional half-perimeter approximations~\cite{hsu2012unified}.
Instead, we define the total wirelength $WL(\mathbf{x}, \mathbf{y}, \mathbf{\theta})$ following its exact definition:
\begin{equation}
WL=\Sigma_{e\in\mathbb{E}}\Sigma_{p_i,p_j\in\mathbb{B}_e}\left(\|X_{p_i}-X_{p_j}\|+\|Y_{p_i}-Y_{p_j}\|\right), \label{equ:wl}
\end{equation}
where $p_i$ and $p_j$ are the two microbumps belonging to net \textit{e} and 
\begin{equation}
    \begin{aligned}
        X_{p_i(j)}&= x_{i(j)}+x_{p_{i(j)}}\cdot\cos\theta_{i(j)}-y_{p_{i(j)}}\cdot\sin\theta_{i(j)},\\
        Y_{p_i(j)}&= y_{i(j)}+x_{p_{i(j)}}\cdot\sin\theta_{i(j)}-y_{p_{i(j)}}\cdot\cos\theta_{i(j)},
    \end{aligned}
\end{equation}
for $i(j)$ representing the corresponding chiplet of microbump $p_i(p_j)$.

To incorporate thermo-mechanical awareness into the optimization framework, we jointly consider thermal distribution and warpage constraints. Let $T_g(\mathbf{x}, \mathbf{y}; \boldsymbol{\theta})$ denote the steady-state temperature field calculated from~\eqref{equ:thermal} parameterized by design variables $\boldsymbol{\theta}$, and warpage be the out-of-plane displacement computed from the thermal profile as defined in~\eqref{equ:warpage}. 
Specifically, we require that: (i) the temperature at \textit{every} grid point inside the discretized grids $\mathcal{G}$ over the interposer surface satisfies $T(\mathbf{r}; \mathbf{x}, \mathbf{y}, \boldsymbol{\theta}) \le T_{\text{th}}$, and (ii) the total package warpage satisfies $\text{Warpage}(\mathbf{x}, \mathbf{y}, \boldsymbol{\theta}) \le W_{\text{th}}$, where $T_{\text{th}}$ and $W_{\text{th}}$ are reliability thresholds, beyond which performance degradation, electromigration, or thermal runaway may occur. 
For multiple workloads, the constraints can be conservatively extended by enforcing the limit among all scenarios—a lightweight generalization that scales linearly with the number of profiles.

To avoid overlapping between chiplets, we employ the bell-shaped density function as in~\cite{hsu2012unified,lin2023routability}. After the placement region (interposer) has been divided into uniform bin grids, the density function $D_b(\mathbf{x}, \mathbf{y}, \mathbf{\theta})$ in bin $b$ is defined as:
\begin{equation}
    D_b(\mathbf{x},\mathbf{y},\theta)=
    \sum_{C_i\in\mathbb{C}}\left(D_i\times P_x(b,C_i) P_y(b,C_i)\right),
    \label{equ:dens}
\end{equation} 
where $D_i$ is the normalization factor, $P_x(b,C_i),\ P_y(b,C_i)$ are the overlap functions of bin b and chiplet $C_i$ along the x and y directions. Since their plain expressions are neither smooth nor differentiable, the bell-shaped potential function is used to smooth $P_x$ and $P_y$. 

With the preparations above, we can model the chiplet placement as a constrained optimization problem:
\begin{equation}
    \begin{array}{cl}
    \min & WL(\mathbf{x}, \mathbf{y}, \mathbf{\theta}) \\
    \text { s.t. } & D_b(\cdot) \leq M_b, \quad \text { for each bin } b, \\
    & T(\mathbf{r}; \cdot) \le T_{\text{th}}, \quad \text { for } r \in \mathcal{G}, \\
    & \text{Warpage}(\cdot) \le W_{\text{th}}
    \end{array}
    \label{equ:constraint}
\end{equation}
where $M_b$ is the maximum allowable area inside $b$, defined as: $M_b = t_{max}(w_bh_b)$, where $t_{max}$ is a parameter representing the target density value for each bin, and $w_b\ (h_b)$ is the width (height) of bin $b$.
Utilizing the quadratic penalty method, it is further translated into a sequence of unconstrained minimization problems:
\begin{equation}
\begin{aligned}
\underset{\mathbf{x,y},\theta}{\min} 
    \{\mathcal{J}(\mathbf{x}, \mathbf{y}, \boldsymbol{\theta})\}, \\
    \mathcal{J}(\mathbf{x}, \mathbf{y}, \boldsymbol{\theta}) = 
    & \underbrace{WL(\mathbf{x}, \mathbf{y}, \boldsymbol{\theta})}_{\text{wirelength}} \\ 
    + & \underbrace{\lambda_{\text{dens}} \sum_b \left( D_b(\mathbf{x}, \mathbf{y}, \boldsymbol{\theta}) - M_b \right)^2}_{\text{density penalty}} \\
    + & \underbrace{\lambda_T \sum_{\mathbf{r} \in \mathcal{G}} \left[ \max\left(0,\, T(\mathbf{r}; \mathbf{x}, \mathbf{y}, \boldsymbol{\theta}) - T_{\text{th}} \right) \right]^\gamma}_{\text{thermal penalty}} \\
    + & \underbrace{\lambda_W \left[ \max\left(0,\, \text{Warpage}(\mathbf{x}, \mathbf{y}, \boldsymbol{\theta}) - W_{\text{th}} \right) \right]^\gamma}_{\text{mechanical penalty}}, \\
    \label{eq:objective}
\end{aligned}
\end{equation}
where $\lambda_T, \lambda_W, \lambda_{\text{dens}} > 0$ are tunable weights balancing wirelength, thermal safety, warpage, and density, respectively; $\gamma \in \mathbb{Z}^+$ controls penalty aggressiveness (e.g., $\gamma=2$ for quadratic, $\gamma>2$ for stronger discouragement of violations). 
When all constraints are satisfied ($T \le T_{\text{th}}$, $\text{Warpage} \le W_{\text{th}}$, $D_b \le M_b$), the penalty terms vanish, and the optimizer minimizes wirelength and congestion—yielding a smooth, gradient-friendly objective that prioritizes feasible, high-performance placements. Otherwise, the penalty actively steers the solution toward feasible thermo-mechanical regimes.
A value of $\gamma > 1$ yields a superlinear penalty growth with increasing temperature or warpage violations, thereby strongly discouraging high-temperature or high-warpage solutions.

Finally, we formulate the problems in this work as follows.

\textbf{Problem I (Compact Model Fitting):}
\textit{Given the physics configuration, including system geometry, material parameters, and power, the goal is to generate a compact thermo-mechanical model that can predict the temperature and displacement distribution as accurately as possible compared to golden simulator}.

\textbf{Problem II (Placement Optimization):}  
\textit{Given the design specifications and chiplet information—including size, microbump locations, and connectivity—the goal is to supports both wirelength-driven and thermo-mechanical-aware objectives, and can efficiently find Pareto-optimal placements that balance wirelength, temperature, and warpage according to the design priority.}

\section{ATMPlace Framework}
\label{sec:Algorithm}

This section presents the \texttt{ATMPlace} framework for thermo-mechanical-aware chiplet placement in 2.5D ICs. The overall architecture is introduced in \secref{sec:framework}, followed by the compact thermal and warpage models in Sections~\ref{sec:compact} and~\ref{sec:warpage}, respectively. The placement algorithm comprises initialization (\secref{sec:initial}), multi-physics-aware optimization (\secref{sec:thermalflow}), and legalization (\secref{sec:legalize}).

\subsection{Framework Overview}
\label{sec:framework}

Figure~\ref{fig:framework} depicts the \texttt{ATMPlace} workflow. First, compact thermal and warpage models are trained offline using high-fidelity simulation data (Sections~\ref{sec:compact}--\ref{sec:warpage}). Subsequently, given a placement instance—defined by chiplet dimensions, microbump connectivity, and interposer specifications—the algorithm proceeds in three stages: (i) an MILP-based initialization generates a non-overlapping, orientation-aware seed layout (\secref{sec:initial}); (ii) a gradient-based optimizer co-minimizes wirelength, thermal hotspots, and warpage using the compact models (\secref{sec:thermalflow}); (iii) a final MILP-based legalization stage eliminates residual overlaps and refines connectivity while preserving physical integrity (\secref{sec:legalize}).

\subsection{Compact Thermal Model}
\label{sec:compact}

Accurate and efficient computation of the temperature gradient $\nabla T$ with respect to chiplet locations is essential for gradient-based placement optimization. However, direct numerical simulation incurs prohibitive cost for large-scale systems~\cite{park1999inverse}. Existing Green’s function-based compact models~\cite{Torregiani2009compact,kung2011thermal} discretize the domain into $M \times M$ grids and compute pairwise responses, yielding $\mathcal{O}(M^4)$ complexity for $N$ chiplets. Moreover, these approaches model power sources as point or uniform-area heat sources and neglect material heterogeneity (e.g., silicon dies, epoxy underfill, Cu pillars), limiting accuracy for large-area, multi-material chiplets.

To overcome these limitations, we propose a chiplet-centric compact thermal model that: (i) treats each chiplet as a spatially extended, piecewise-homogeneous heat source; (ii) accounts for interfacial thermal resistance via effective length normalization; and (iii) reduces complexity to $\mathcal{O}(N M^2)$, where $N \ll M^2$ in practice.

We begin with the steady-state heat equation under quasi-homogeneous approximation:
\begin{equation}
    \nabla^2 T(\mathbf{r}) = - \frac{P(\mathbf{r})}{\kappa(\mathbf{r})},
    \label{eq:thermal_poisson}
\end{equation}
where $T(\mathbf{r})$ is temperature, $P(\mathbf{r})$ is volumetric power density, and $\kappa(\mathbf{r})$ is thermal conductivity. Let $G(\mathbf{r}, \mathbf{r}')$ denote the free-space Green’s function satisfying $\nabla^2 G(\mathbf{r}, \mathbf{r}') = \delta(\mathbf{r} - \mathbf{r}')$, with solution $G(\mathbf{r}, \mathbf{r}') = \frac{1}{4\pi \|\mathbf{r} - \mathbf{r}'\|}$. Applying Green’s theorem yields the integral formulation:
\begin{equation}
    T(\mathbf{r}) = -\int_{\Omega} G(\mathbf{r}, \mathbf{r}') \frac{P(\mathbf{r}')}{\kappa(\mathbf{r}')} \, d\mathbf{r}'
    = \sum_{i=1}^N \frac{P_i}{4\pi \kappa_i} \int_{\Omega_i} \frac{d\mathbf{r}'}{\|\mathbf{r} - \mathbf{r}'\|},
\end{equation}
where $\Omega_i$ is the 3D domain of chiplet $i$ with uniform power density $P_i$ and conductivity $\kappa_i$. The free-space assumption is justified by the interposer’s lateral dimensions being significantly larger than individual chiplet footprints.
To evaluate the volume integral efficiently, we apply the identity~\cite{ciftja2020electrostatic}:
\begin{align}
    & \frac{1}{\|\mathbf{r} - \mathbf{r}'\|} 
    = \frac{2}{\sqrt{\pi}} \int_0^\infty e^{-u^2 \|\mathbf{r} - \mathbf{r}'\|^2} \, du 
    \quad \notag \\
    = & \frac{2}{\sqrt{\pi}} \int_0^\infty \!\!\iiint_{\Omega_i} 
       e^{-u^2 \left[(x - x')^2 + (y - y')^2 + (z - z')^2\right]} 
       dx'\,dy'\,dz' \, du. 
    \label{eq:u_integral}
\end{align}

Assuming chiplet thickness $t_i$ is uniform and $t_i \ll w_i, h_i$, we approximate the $z'$-integral as $t_i e^{-u^2 a^2}$, where $a \in [0, t_i]$ is an effective decay depth (fitted during training). Utilizing the error function $\textbf{erf}\left(x\right)=\frac{2}{\sqrt{\pi}}\int_0^xe^{-t^2}dt$, the $x'$ and $y'$ integrals admit closed-form solutions $K(u,x,x_i,w_i)\triangleq$:
\begin{equation}
    \frac{\sqrt{\pi}}{2u}[\textbf{erf}\left(u(\frac{w_i^\prime}{2}-(x-x_i))\right)+\textbf{erf}\left(u(\frac{w_i^\prime}{2}+(x-x_i))\right)].
\end{equation}
and similarly for $y'$. Substituting into~\eqref{eq:u_integral} yields:
\begin{equation}
    \frac{\sqrt{\pi} t_i}{2} \int_0^\infty \frac{e^{-u^2 a^2}}{u^2} \, K(u, x, x_i, w_i) K(u, y, y_i, h_i) \, du.
\end{equation}

This integral reduces to a linear combination of four instances of the auxiliary function~\cite{ciftja2020electrostatic}:
\begin{equation}
    F(a, b, c) = \int_0^\infty \frac{e^{-a^2 x^2}}{x^2} \operatorname{erf}(b x) \operatorname{erf}(c x) \, dx,
\end{equation}
with closed-form expression:
\begin{equation}
    F(a, b, c) = \frac{2}{\sqrt{\pi}} \Bigg[
        b \ln \frac{c + \Delta}{\sqrt{a^2 + b^2}} +
        c \ln \frac{b + \Delta}{\sqrt{a^2 + c^2}} -
        a \tan^{-1} \frac{bc}{a \Delta}
    \Bigg],
    \label{eq:F_closed}
\end{equation}
where $\Delta = \sqrt{a^2 + b^2 + c^2}$. To model heterogeneous thermal interfaces, we introduce anisotropic length normalization $l_{x,i}, l_{y,i} > 0$ per chiplet, yielding the final compact thermal model:
\begin{align}
    T_c(x, y) = \sum_{i=1}^N A P_i \bigg[ &
    F\!\bigg(a,\, \tfrac{w_i/2 - (x - x_i)}{l_{x,i}},\, \tfrac{h_i/2 - (y - y_i)}{l_{y,i}} \bigg) \nonumber \\
    + & F\!\bigg(a,\, \tfrac{w_i/2 - (x - x_i)}{l_{x,i}},\, \tfrac{h_i/2 + (y - y_i)}{l_{y,i}} \bigg) \nonumber \\
    + & F\!\bigg(a,\, \tfrac{w_i/2 + (x - x_i)}{l_{x,i}},\, \tfrac{h_i/2 - (y - y_i)}{l_{y,i}} \bigg) \nonumber \\
    + & F\!\bigg(a,\, \tfrac{w_i/2 + (x - x_i)}{l_{x,i}},\, \tfrac{h_i/2 + (y - y_i)}{l_{y,i}} \bigg) \bigg] + B,
    \label{eq:compact_thermal}
\end{align}
where $A$ (amplitude), $a$ (decay depth), and $B$ (bias) are global parameters. Thus, the model uses $2N + 3$ trainable parameters.

Parameter fitting is performed offline: given a set of legal placements (typically 5–15, sufficient for convergence), ground-truth temperature fields $\{T_{\text{label}}\}$ are generated using numerical solver. The parameters $\boldsymbol{\beta} = \{A, a, B, \{l_{x,i}, l_{y,i}\}_{i=1}^N\}$ are optimized via: $\arg\min_{\boldsymbol{\beta}} \big\| T_c(\boldsymbol{\beta}) - T_{\text{label}} \big\|_2^2$
implemented in PyTorch. The closed-form structure of~\eqref{eq:compact_thermal} ensures rapid convergence and robustness to initialization.

\subsection{Compact Warpage Model}
\label{sec:warpage}
Full FEA-based warpage prediction is computationally infeasible for iterative placement. Inspired by the beam-theory analysis of~\citet{tasi2020}, we derive a compact, differentiable surrogate model.
Consider a three-layer package stack. The vertical displacement $D^z(r)$ along radial coordinate $r$ (origin at package centroid) is:
\begin{align}
    D_1^z(r) &= \frac{t \Delta \alpha \Delta T}{2 \lambda D} \left( \tfrac{1}{2} r^2 - \tfrac{\cosh kr - 1}{k^2 \cosh kl_1} \right), \quad 0 \leq r \leq l_1, \\
    D_2^z(r) &= D_1^z(l_1) + D_1^{z'}(l_1)(r - l_1) - \tfrac{1}{2} \nu_3 D_1^{z''}(0) (r - l_1)^2 \notag \\
    &\quad l_1 \leq r \leq l_2,
\end{align}
where $t$ is total thickness, $\Delta \alpha$ is CTE difference, $\Delta T$ is thermal excursion, $\nu_i$ and $E_i$ are Poisson’s ratio and Young’s modulus, $D_i = E_i t_i^3 / [12(1 - \nu_i)]$ is flexural rigidity, and $\lambda = \frac{1 - \nu_1}{E_1 t_1} + \frac{1 - \nu_3}{E_3 t_3} + \frac{t^2}{4D}$. Empirical studies indicate the second-order term in $D_2^z$ contributes less than 5\% to total warpage; thus, we approximate the displacement field as a localized quadratic form.

For each chiplet $C_i$, the warpage is modeled as a spatially varying quadratic–hyperbolic field centered at $(x_i, y_i)$, where the anisotropy and effective decay in the $x$- and $y$-directions are captured by learnable scaling factors. Inspired by the analytical form in \citet{tasi2020}, but adapted for 2D in-plane placement and computational efficiency, we define the \textit{local} warpage contribution of chiplet $i$ as:
\begin{equation}
\begin{aligned}
    w_i(x, y)
    &= \Bigl[\, k_{x,i}^2 (x - x_i)^2 + k_{y,i}^2 (y - y_i)^2 \Bigr] \\
    & + \lambda_i \Bigl[\, k_{x,i} (x - x_i) + k_{y,i} (y - y_i) \Bigr] 
      + c_i,
\end{aligned} \label{eq:w_local}
\end{equation}
where $k_{x,i}, k_{y,i} > 0$ are anisotropic decay/curvature parameters (learnable), controlling the spatial \textit{spread} of the warpage field—effectively normalizing physical distances by material-dependent thermal diffusion or mechanical stiffness scales; $\lambda_i$ is a linear coupling coefficient, allowing asymmetry in the warpage gradient (e.g., due to boundary constraints or non-uniform underfill curing); and $c_i$ is a constant offset representing baseline deformation associated with chiplet $i$.

The total warpage field is then obtained by thermally weighting each local contribution and summing:
\begin{equation}
W(x, y) = \alpha \cdot \sum_{i=1}^{N} 
\bigl(T(x, y) - T_{\mathrm{ref}}\bigr) \cdot w_i(x, y) \;+\; b
\label{eq:warpage_compact}
\end{equation}
Here, $\alpha$ is a global amplitude factor (shared across chiplets), capturing the overall thermo-mechanical coupling strength; $T_{\mathrm{ref}}$ is the reference temperature, accounting for process-induced pre-stress or ambient thermal history; and $b$ is a global bias term, representing uniform out-of-plane offset (e.g., due to global CTE mismatch or post-molding shrinkage). Critically, the term $(T(x, y) - T_{\mathrm{ref},\,i})$ takes the runtime local heating effects into account.
Crucially, $W(x, y)$ is fully differentiable with respect to chiplet locations $(x_i, y_i)$ and orientations (via $T_c$ and $w_i$).

Model parameters $\boldsymbol{\gamma} = \{\alpha, b, \{k_{x,i}, k_{y,i}, \lambda_i, c_i, T_{\text{ref},i}\}_{i=1}^N\}$ (total $5N + 2$) are fitted against FEA-generated warpage maps $\{W_{\text{label}}\}$: $\arg\min_{\boldsymbol{\gamma}} \big\| W(\boldsymbol{\gamma}) - W_{\text{label}} \big\|_2^2$.
Convergence is accelerated by the physical priors embedded in~\eqref{eq:w_local}--\eqref{eq:warpage_compact}:
\textbf{Localized curvature} around each chiplet (via the quadratic term), proportional to local temperature rise; \textbf{Asymmetric tilting or shear} (via the linear term), arising from heterogeneous boundary conditions or placement-dependent stress concentration with only $5N + 3$ learnable parameters (including $\alpha$ and $b$).

\subsection{Initialization}\label{sec:initial}
A high-quality initial placement—jointly optimizing chiplet positions and orientations—is critical for convergence and solution quality~\cite{agnesina2023autodmp}. To address this, we propose a MILP-based initializer that explicitly incorporates rotation and microbump clustering.

Let $\mathbb{A}$ denote the set of connected chiplet pairs. For each pair $(C_i, C_j) \in \mathbb{A}$, let $A_{ij}$ be the number of D2D nets between them. Due to the structured layout of interface bumps (e.g., UCIe), these connections form spatially compact ``clumps'' rather than being uniformly distributed. We precompute the centroid offsets of each clump relative to its chiplet center: $(O_{ij}^x, O_{ij}^y)$ on $C_i$ and $(O_{ji}^x, O_{ji}^y)$ on $C_j$, directly from the bump map.
Each chiplet $C_i$ supports four legal orientations: $\theta_i \in \{0^\circ, 90^\circ, 180^\circ, 270^\circ\}$. We encode $\theta_i$ using two binary variables $u_i, v_i \in \{0,1\}$ via the bijection:
\begin{equation*}
    \{0^\circ, 90^\circ, 180^\circ, 270^\circ\}
    \leftrightarrow 
    \{(0,0), (0,1), (1,1), (1,0)\}.
\end{equation*}
The absolute position of the clump on $C_i$ is then:
\begin{align*}
    X_{ij} &= x_i + (1 - u_i - v_i) O_{ij}^x - (v_i - u_i) O_{ij}^y, \\
    Y_{ij} &= y_i + (v_i - u_i) O_{ij}^x + (1 - u_i - v_i) O_{ij}^y,
\end{align*}
which yields linear expressions for all rotations—enabling exact MILP embedding.
Using these, we define a bump-aware wirelength proxy for initialization:
\begin{equation}
    \text{WL}^{(\text{init})} = \sum_{(i,j) \in \mathbb{A}} A_{ij} \left( |X_{ij} - X_{ji}| + |Y_{ij} - Y_{ji}| \right),
    \label{eq:init_wl}
\end{equation}
where the Manhattan distance approximates the total routed length of all $A_{ij}$ nets between the two clumps.

Subject to the following constraints:  
(i) \textit{Boundary containment}:  
\begin{equation}
    \frac{w_i'}{2} \leq x_i \leq W - \frac{w_i'}{2}, \quad
    \frac{h_i'}{2} \leq y_i \leq H - \frac{h_i'}{2},
    \label{eq:bound}
\end{equation}
with rotated dimensions $w_i' = |1 - u_i - v_i|\, w_i + |v_i - u_i|\, h_i$, $h_i' = |v_i - u_i|\, w_i + |1 - u_i - v_i|\, h_i$.  
(ii) \textit{Non-overlap}: enforced via four disjunctive constraints per chiplet pair $(i,j)$ using auxiliary binaries $\delta_{ij}^{(k)} \in \{0,1\}$ ($k = 1,\dots,4$):
\begin{equation}
    \begin{aligned}
        x_i + (w'_i + w'_j)\varepsilon &\le x_j + W \delta^{(1)}_{ij}, \\
        x_j + (w'_i + w'_j)\varepsilon &\le x_i + W \delta^{(2)}_{ij}, \\
        y_i + (h'_i + h'_j)\varepsilon &\le y_j + H \delta^{(3)}_{ij}, \\
        y_j + (h'_i + h'_j)\varepsilon &\le y_i + H \delta^{(4)}_{ij},
    \end{aligned}\label{eq:overlap}
\end{equation}
where $\varepsilon \in [0, 0.5)$ is a small slack factor to ensure strict separation (e.g., $\varepsilon = 0.1$ yields $10\%$ bin-width gap). Feasibility is guaranteed by:
\begin{equation}
    \sum_{k=1}^4 \delta_{ij}^{(k)} \le 3.
    \label{eq:delta}
\end{equation}

The full initialization problem is:
\[
\min_{\mathbf{x}, \mathbf{y}, \mathbf{u}, \mathbf{v}, \boldsymbol{\delta}} \text{WL}^{(\text{init})} \quad \text{s.t.} \quad \eqref{eq:bound}\text{--}\eqref{eq:delta} \text{\quad holds},
\]
solved efficiently using an off-the-shelf MILP solver. Though MILP scales exponentially in theory, our formulation remains tractable for dozens of chiplets—serving as a robust warm-start for subsequent gradient-based refinement.

\subsection{Multi-Physics-Aware Optimization}
\label{sec:thermalflow}

Starting from the initialized layout, we perform gradient-based co-optimization of wirelength, density, temperature, and warpage. The objective is:
\begin{equation}
    \mathcal{F}(\mathbf{X}) = \text{WL}' + \lambda_{\text{dens}} \mathcal{D} + 
    \lambda_T\|T_c(\cdot)-T_{\text{th}}\|_2^\gamma
    +\lambda_W\|W_c(\cdot)-W_{\text{th}}\|_2^\gamma,
\end{equation}
where $\mathbf{X} = \{(x_i, y_i, \theta_i)\}_{i=1}^N$, $\text{WL}'$ is the orientation-aware wirelength (see below), $\mathcal{D} = \sum_b (D_b' - M_b)^2$ penalizes bin density deviation from target $M_b$, $T_c(\cdot)$ and $W_c(\cdot)$ are the outputs of the compact thermal and warpage models, and they are differentiable w.r.t. $x_i,y_i$ and $\theta_i$ through the chain rule. 

\subsubsection{Orientation-aware wirelength and density}
Our orientation-aware optimization algorithm is based on the model proposed by \citet{lin2023routability}. 
Chiplets admit discrete rotations $\Theta = \{0^\circ, 90^\circ, 180^\circ, 270^\circ\}$. We use a smooth projection that calculates the probability of the chiplet $C_i$ to be rotated to each legal orientation $\theta^k$~\cite{lin2023routability}:
\begin{equation}
    B_z(\theta^k, \theta_i) = 
    \frac{\exp\big(R_z(\theta^k, \theta_i)/\eta\big)}{\sum_{\theta^k \in \Theta} \exp\big(R_z(\theta^k, \theta_i)/\eta\big)},
\end{equation}
where $\eta > 0$ controls sharpness. 
With angular deviation $\Delta\theta_i^k$, $R_z(\theta^k, \theta_i)$ defined piecewise as:
\begin{equation}
    \begin{aligned}
    \Delta\theta_i^k=&\left|0.5-\left|0.5-\|\theta_i-\theta^k\|/{360}\right|\right|,\\
    R_z\left(\theta^k, \theta_i\right)=& 
    \begin{cases}
        1-2|\Theta|^2\left|\Delta\theta_i^k\right|^2, 
        & 0<\Delta\theta_i^k \leq \frac{1}{2|\Theta|} \\ 
        2|\Theta|^2\left(\Delta\theta_i^k-\frac{1}{|\Theta|}\right)^2, 
        & \frac{1}{2|\Theta|}<\Delta\theta_i^k \leq \frac{1}{|\Theta|}\\
        0, & \text { otherwise, }
    \end{cases}\\
    \end{aligned}
\end{equation}

The projected wirelength and density become:
\begin{align}
    \text{WL}' &= \sum_{e \in \mathbb{E}} \sum_{p_i, p_j \in \mathcal{B}_e} \sum_{\theta^k, \theta^\ell \in \Theta} \!\! B_z(\theta^k, \theta_i) B_z(\theta^\ell, \theta_j) \, \text{HPWL}(p_i^{(k)}, p_j^{(\ell)}), \\
    D_b' &= \sum_{C_i \in \mathbb{C}} \sum_{\theta^k \in \Theta} B_z(\theta^k, \theta_i) \cdot \text{Area}(C_i^{(k)} \cap \text{bin } b).
\end{align}

\subsubsection{Optimization algorithm}
In our framework, the minimization of \eqref{eq:objective} is solved by the conjugate gradient descent (CGD) method. 
The algorithm is summarized in Algorithm \ref{alg:cgd} with adaptive step sizing and noise injection to escape local minima. 
Step sizes for position and angle are decoupled: $\alpha_k = \big(\alpha_{k,x}, \alpha_{k,y}, \alpha_{k,\theta}\big)$. The density weight $\lambda_{\text{dens}}$ is adaptively scaled:
\begin{equation}
    \lambda_{\text{dens}}^{(0)} = \frac{\sum |\nabla_{\mathbf{x},\mathbf{y}} \text{WL}'|}{\sum |\nabla_{\mathbf{x},\mathbf{y}} \mathcal{D}|}, \quad
    \lambda_{\text{dens}}^{(k+1)} = \lambda_{\text{dens}}^{(k)} \cdot \left(1 + \rho \cdot \text{OVFL}^{(k)}\right).
\end{equation}
Finally, the system is prone to get stuck in local optima during the nonlinear optimization. Inspired by \cite{Xue2024EscapingLO}, we perturb the placement by adding random noise when the overflow of the system is stuck at high values or converged. The overflow metric reads:
\begin{equation}
    \text{OVFL} = \frac{1}{S} \sum_b \max\big(0, D_b' - M_b \big),
    \label{eq:overflow}
\end{equation}
where $S$ is the total chiplet area.

\begin{algorithm}[H]
    \caption{Optimization via CGD}\label{alg:cgd}
    \begin{algorithmic}[1]
    \Require $\mathcal{F}(\textbf{X})$: objective function; 
    $Max\_iter$: number of max iteration; $\textbf{lr}=(lr_{pos},lr_{ang})$: learning rate of position and angular variables
    \Ensure optimal $\textbf{X}^*$
    \State initialize $\textbf{X}_0$ (\secref{sec:initial}), $\lambda_{dens}$, $\textbf{g}_0=\textbf{0}$, and $\textbf{d}_0=\textbf{0}$;
    \For{$k = 1$ to $Max\_iter$}
        \State Compute $\mathcal{F}(\mathbf{X}_{k-1})$ and gradient $\mathbf{g}_k \gets \nabla \mathcal{F}(\mathbf{X}_{k-1})$
        \State $\beta_k \gets \frac{\mathbf{g}_k^\top (\mathbf{g}_k - \mathbf{g}_{k-1})}{\|\mathbf{g}_{k-1}\|_2^2 + \epsilon}$ \quad (Polak–Ribière, $\beta_1 = 0$)
        \State $\mathbf{d}_k \gets -\mathbf{g}_k + \beta_k \mathbf{d}_{k-1}$
        \State $\alpha_k \gets \big( \eta_x / \|\mathbf{d}_k^{(x)}\|, \, \eta_y / \|\mathbf{d}_k^{(y)}\|, \, \eta_\theta / \|\mathbf{d}_k^{(\theta)}\| \big)$
        \State $\mathbf{X}_k \gets \mathbf{X}_{k-1} + \alpha_k \odot \mathbf{d}_k$
    \EndFor
    \State legalize (\secref{sec:legalize}) and derive the final result $\textbf{X}^*$
    \end{algorithmic}
\end{algorithm}

\subsection{Legalization}\label{sec:legalize}
After optimization, the chiplets are roughly uniformly distributed across the interposer. But there are still some overlapping areas at this time. Additionally, sometimes minor temperature optimization can lead to excessive increases in total wirelength. 
Thus, a legalization stage is indispensable to eliminate all overlaps and reduce the total wirelength while keeping the changes in chiplet positions minimal, essentially maintaining the temperature distribution. 
In this stage, given the optimized solution $\left(x_i^{(opt)}, y_i^{(opt)}, \theta_i^{(opt)}\right)$, we fix the orientations and optimize the positions of chiplets by MILP.

The optimization objective includes two parts. The first term seeks to minimize the total wirelength:
\begin{align}
    WL^{(legal)}=
    \Sigma_{e\in\mathbb{E}}\Sigma_{p_i,p_j\in\mathbb{B}_e}
    \bigg(&\|(x_i+x_{p_i})-(x_j+x_{p_j})\|\notag \\
    +&\|(y_i+y_{p_i})-(y_j+y_{p_j})\|\bigg), 
\end{align}
While the second term constrains the displacement:
\begin{equation}
    DSP=
    \Sigma_{C_i}
    \left(\|x_i-x_i^{(opt)}\|+\|y_i-y_i^{(opt)}\|\right).
\end{equation}
Combined with the non-overlapping conditions formulated before (\eqref{eq:overlap} with $\varepsilon=1$) with the chiplet width and height modified to involve the minimum spacing between chiplets, the final legalization problem is:
\[
\min_{\mathbf{x}, \mathbf{y}} DSP+\lambda_w\cdot\ WL^{(legal)} \quad \text{s.t.} \quad \eqref{eq:bound}\text{--}\eqref{eq:delta} \text{\quad holds},
\]
Here $\lambda_w$ is a parameter that controls the extent to which the total wirelength is optimized in this stage. Since it may be time-consuming or sometimes even infeasible to optimize the wirelength, the algorithm will set $\lambda_w$ to $0$ when the time limit of 10 minutes reached.

Ending on a note, experiments have shown that the quality of the placement results heavily relies on the choice of parameters. Thus, inspired by \texttt{AutoDMP} \cite{agnesina2023autodmp}, which puts forth to optimize the hyperparameters through Bayesian optimization to improve overall performance, we also carry out hyperparameter optimization to find the most favorable parameters. 
\section{Experimental Results}
\label{sec:Results}

This section presents a comprehensive evaluation of the proposed \texttt{ATMPlace} framework. We first describe the experimental setup in~\secref{sec:setup}, followed by validation of the compact thermal and warpage models. A detailed analysis of placement quality, runtime, and multi-objective trade-offs is provided in~\secref{sec:placement}.

\subsection{Experimental Setup}
\label{sec:setup}

Given the scarcity of publicly available large-scale 2.5D-IC benchmarks and the omission of die-to-die (D2D) interconnects in prior works, we construct ten representative test cases, summarized in Table~\ref{tab:cases}. These cases span 6 to 61 chiplets, up to thousands of nets, and whitespace ratios between 0.35 and 0.65. The first five cases are high-performance computing (HPC) systems comprising CPUs, GPUs, HBMs, and DRAMs, adapted from~\cite{ma2021tap}. The latter five incorporate analog and micro-electromechanical system (MEMS) modules~\cite{10069926,nasrullah2019designing} to evaluate heterogeneous integration capability. Chiplet power densities range from $\SI{2e5}{\watt\per\square\meter}$ to $\SI{3e6}{\watt\per\square\meter}$.

For thermal simulation, we use \textbf{ATSim}~\cite{wang2024atsim} as the golden reference due to its capability in modeling multiple thermal properties, like anisotropy and heterogeneity. 
Mechanical warpage is simulated using an in-house finite-element solver~\cite{zhu2025warpage}. Both simulators are verified and calibrated against the ANSYS commercial simulators.
Material parameters and boundary conditions follow established conventions in the literature~\cite{11101154,zhu2025warpage}. 
A uniform grid resolution of $64 \times 64$ is employed.

Inter-chiplet communication is enabled via D2D link modules. We adopt both the \textit{standard} ($\times16$) and \textit{advanced} ($\times32$) UCIe-compliant interfaces.
The $\times16$ module offers longer channel reach, while the $\times32$ variant provides higher bandwidth. Key interface parameters are listed in Table~\ref{tab:interface}, with bump pitch geometries illustrated in Fig.~\ref{fig:bump}.

The \texttt{ATMPlace} framework\footnote{Benchmarks and binaries are open-sourced at \url{https://github.com/Brilight/ATPlace_pub}.} is implemented in Python, building upon \texttt{DREAMPlace}~\cite{lin2019dreamplace} for infrastructure and \texttt{Optuna}~\cite{optuna} for hyperparameter optimization. MILP subproblems are solved using \texttt{Gurobi}~\cite{gurobi}. 

All experiments are conducted on a Linux server equipped with dual Intel Xeon Silver 4210 CPUs (2.20~GHz, 20 cores total), 128~GB RAM, and a maximum time budget of 12~hours per run. For fairness, all baseline methods—including enumeration-based \texttt{SP-CP}—are allowed up to 80 CPU cores and the same time limit.

\begin{table}[tbh]
    \centering
    \caption{Benchmark configurations.}
    \label{tab:cases}
    \normalsize
    \resizebox{0.48\textwidth}{!}{
    \begin{tabular}{|c|c|cc|ccc|}
    \hline
    \multirow{2}{*}{\textbf{Case}} 
    & \multirow{2}{*}{\textbf{Bump Type}}
    & \multirow{2}{*}{\textbf{Dies}}
    & \multirow{2}{*}{\textbf{Nets}}
    & \multicolumn{3}{c|}{\textbf{Interposer}}  \\
     & & & & Width (\SI{}{\mm}) & Height (\SI{}{\mm}) & Whitespace/\% \\ \hline
    \textbf{1} & $\times32$ & 6  & 3168 & 42.0 & 42.0 & 40 \\ 
    \textbf{2} & $\times32$ & 6  & 3520 & 55.0 & 52.0 & 65 \\ 
    \textbf{3} & $\times32$ & 8  & 8448 & 39.0 & 39.0 & 60 \\ 
    \textbf{4} & $\times32$ & 11 & 7040 & 57.0 & 59.0 & 40 \\ 
    \textbf{5} & $\times32$ & 12 & 7392 & 37.0 & 37.0 & 35 \\ 
    \textbf{6} & $\times16$ & 20 & 5632 & 49.0 & 53.0 & 55 \\
    \textbf{7} & $\times16$ & 28 & 2816 & 30.0 & 25.0 & 55 \\
    \textbf{8} & $\times16$ & 36 & 2948 & 26.0 & 23.0 & 60 \\ 
    \textbf{9} & $\times16$ & 44 & 7656 & 59.0 & 55.0 & 45 \\     
    \textbf{10}& $\times16$ & 61 & 5280 & 47.0 & 47.0 & 50 \\  
    \hline
    \end{tabular}
    }
\end{table}

\begin{table}[tbh]
    \centering
    \caption{Characteristics of UCIe-compliant interfaces.}
    \label{tab:interface}
    \normalsize
    \resizebox{0.48\textwidth}{!}{
    \begin{tabular}{|c|cc|ccc|}
    \hline
        Package & Cols & Lanes & Bump Pitch ($\SI{}{\mu m}$) & Pitch$_x$ (\SI{}{\mu m}) & Pitch$_y$ (\SI{}{\mu m}) \\ \hline
        Standard\ ($\times16$) & 12 & 16 & 100 & 180 & 90 \\
        Advanced\ ($\times32$) & 16 & 32 & 25 & 27 & 42 \\  \hline
    \end{tabular}
    }
\end{table}

\subsection{Compact Models}

We first assess the fidelity of the proposed compact models, which serve as the foundation for gradient-based placement optimization. For each benchmark case, we generate 20 random placements for training and testing. Model accuracy is quantified via mean absolute error (MAE) and Pearson correlation coefficient ($\rho$). While MAE measures absolute deviation, $\rho$ evaluates not only global physical field alignment but also gradient similarity—critical for optimization stability and convergence.

As shown in Table~\ref{tab:compact}, the compact thermal model achieves an average MAE of $\SI{2.4}{\degreeCelsius}$ and $\rho = 0.977$, indicating excellent agreement with \texttt{ATSim}. Compared to the golden simulator (average runtime $\SI{66}{\second}$ per layout), our model reduces inference time to $\SI{7.5}{\milli\second}$ on average—a speedup exceeding $8000\times$. Similarly, the compact warpage model attains an average MAE of $\SI{2.5}{\micro\meter}$ and $\rho = 0.978$, with a $>10000\times$ speedup over the numerical mechanical simulator.

\begin{figure}[hbt]
\centering
    \begin{subfigure}[t]{0.95\linewidth}
        \centering
        \includegraphics[width=\linewidth]{figs/temp-warp-label.pdf}
        \caption{Golden truth}
        \label{fig:thermal-gt}
    \end{subfigure}
    \begin{subfigure}[t]{0.95\linewidth}
        \centering
        \includegraphics[width=\linewidth]{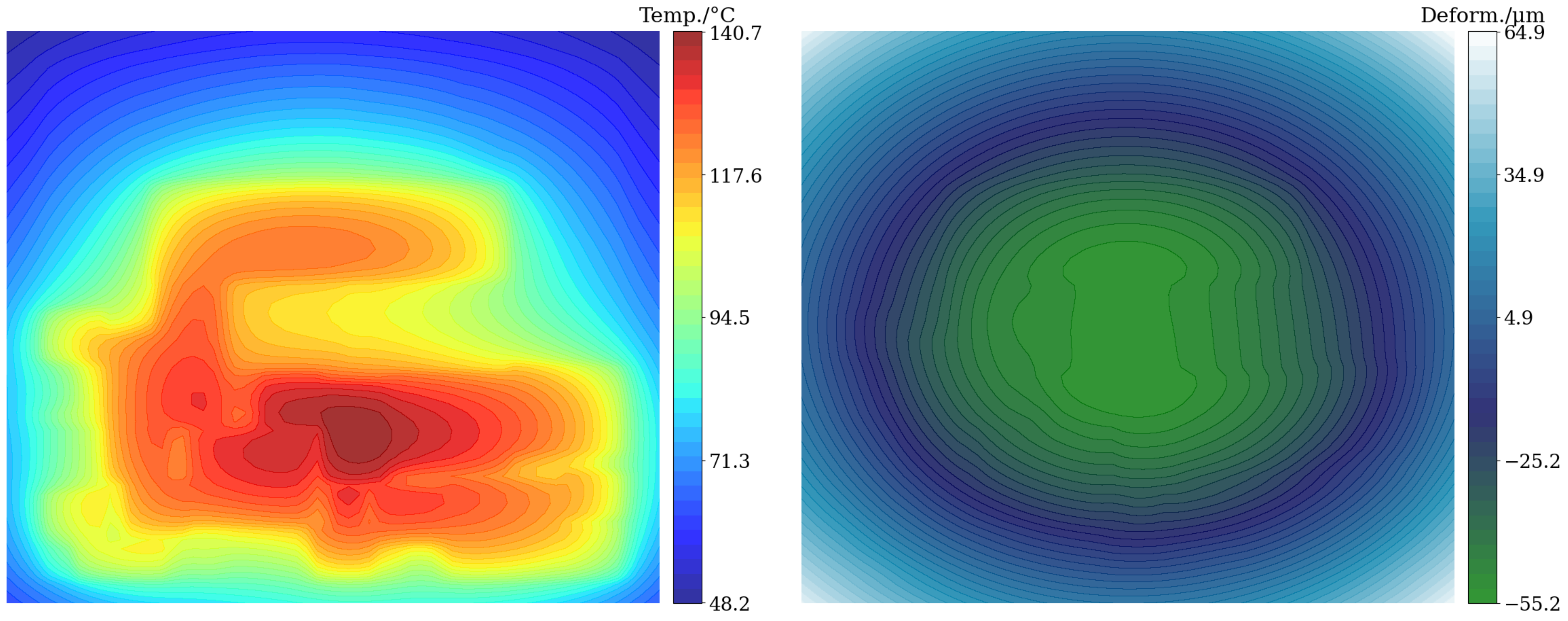}
        \caption{Predictions}
        \label{fig:thermal-pred}
    \end{subfigure}
\caption{Thermo-mechanical field comparison for a placement in Case~1: (a) ground truth, (b) compact model prediction.}
\label{fig:thermal}
\end{figure}

Figure~\ref{fig:thermal} visualizes the temperature and displacement fields for Case~1. The predicted temperature distributions closely match the ground truth, with errors primarily localized at: (1) chiplet edges—attributable to the quasi-homogeneous approximation in ~\eqref{eq:thermal_poisson}; and (2) interposer boundaries—stemming from the infinite-space Green’s function assumption. The latter could be mitigated via image-source methods (e.g., \cite{lin2018fast}), albeit at increased computational cost.

For the displacement field, the dominant source of prediction error arises from the modeling simplification of the linear dependence on local temperature fields~\eqref{eq:warpage_compact}, which neglects higher-order thermo-mechanical coupling. In reality, the displacement field tends to be spatially smoother than the raw temperature field due to structural continuity—our current linear surrogate, while efficient, cannot fully capture this regularization effect.

\begin{table}[tbh]
\centering
\caption{Performance of the proposed compact models. RT denotes runtime; Corr. denotes the correlation coefficient.}
\label{tab:compact}
\normalsize
\resizebox{0.49\textwidth}{!}{
\begin{tabular}{|c|c|ccc|c|ccc|}
\hline
     &\multicolumn{4}{c|}{Thermal} & \multicolumn{4}{c|}{Mechanical} \\ \hline
    \multirow{2}{*}{\textbf{Case}} & \texttt{Golden} & \multicolumn{3}{c|}{Ours} & \texttt{Golden} & \multicolumn{3}{c|}{Ours} \\
    & RT/s & RT/ms & MAE/$\SI{}{\degreeCelsius}$ & Corr. & RT/s & RT/ms & MAE/$\SI{}{\um}$ & Corr. \\ \hline
    \textbf{1} & 70 & 4.7 & 2.15 & 0.984 & 5 & 0.3 & 2.03 & 0.989 \\
    \textbf{2} & 67 & 4.9 & 0.97 & 0.991 & 7 & 0.3 & 2.70 & 0.979 \\
    \textbf{3} & 82 & 6.7 & 3.13 & 0.984 & 5 & 1.1 & 2.17 & 0.973 \\
    \textbf{4} & 78 & 5.9 & 3.18 & 0.977 & 8 & 0.5 & 3.16 & 0.989 \\
    \textbf{5} & 87 & 6.3 & 3.95 & 0.982 & 5 & 0.5 & 2.27 & 0.980 \\
    \textbf{6} & 53 & 10.2 & 2.70 & 0.965 & 8 & 0.7 & 2.77 & 0.979 \\
    \textbf{7} & 58 & 11.4 & 1.92 & 0.955 & 6 & 0.6 & 1.19 & 0.969 \\
    \textbf{8} & 56 & 19.1 & 2.18 & 0.960 & 6 & 1.2 & 0.65 & 0.986 \\
    \textbf{9} & 49 & 16.0 & 2.02 & 0.992 & 13 & 0.7 & 4.07 & 0.985 \\
    \textbf{10} & 60 & 22.1 & 2.13 & 0.981 & 12 & 0.8 & 4.14 & 0.954 \\
    \hline
    \textbf{Avg.} & \textbf{8680}$\times$ & $\textbf{1}\times$ & \textbf{2.43} & \textbf{0.977} & \textbf{13050}$\times$ & $\textbf{1}\times$  & \textbf{2.52} & \textbf{0.978} \\ 
    \hline
\end{tabular}
}
\end{table}

\begin{table*}[htb]
    \centering
    \caption{Comparison of placement algorithms under wirelength-driven optimization. 
    `Wpg' represents the warpage. Normalized metrics use \texttt{ATMPlace} as baseline (1$\times$).}
    \label{tab:wlopt}
    \normalsize
    \resizebox{0.99\textwidth}{!}{
    \begin{tabular}{|c|cccc|cccc|cccc|cccc|}
    \hline
    \multirow{2}{*}{\textbf{Case}} 
    & \multicolumn{4}{c|}{\texttt{TAP-2.5D}~\cite{ma2021tap}}
    & \multicolumn{4}{c|}{\texttt{TACPlace}~\cite{yu2025tacplace}}
    & \multicolumn{4}{c|}{\texttt{SP-CP}~\cite{chiou2023chiplet}}
    & \multicolumn{4}{c|}{\texttt{ATMPlace (WL-driven)}} \\
    &RT/min &TWL/\SI{}{m} &$\text{T}_{\text{max}}/\SI{}{\degreeCelsius}$ &$\text{Wpg}/\SI{}{\um}$
    &RT/min &TWL/\SI{}{m} &$\text{T}_{\text{max}}/\SI{}{\degreeCelsius}$ &$\text{Wpg}/\SI{}{\um}$
    &RT/min &TWL/\SI{}{m} &$\text{T}_{\text{max}}/\SI{}{\degreeCelsius}$ &$\text{Wpg}/\SI{}{\um}$
    &RT/min &TWL/\SI{}{m} &$\text{T}_{\text{max}}/\SI{}{\degreeCelsius}$ &$\text{Wpg}/\SI{}{\um}$ \\ \hline
    \textbf{1} 
    & 3.1 & 27.499 & 101.5 & 84.6
    & 2.0 & 60.537 & 100.2 & 86.4
    & 0.02 & \textbf{11.902} & 101.1 & 84.9 
    & 0.3  & 11.973 & 102.7 & 85.1 \\
    \textbf{2} 
    & 3.5 & 62.124 & 66.4 & 97.4 
    & 2.1 & 74.656 & 66.0 & 97.6
    & 0.01 & 16.107 & 67.4 & 103.8
    & 0.3 & \textbf{15.408} & 67.8 & 105.3 \\
    \textbf{3} 
    & 7.3 & 87.543 & 132.8 & 63.7 
    & 5.5 & 122.972 & 126.1 & 60.9
    & 0.5 & \textbf{32.200} & 144.2 & 66.7 
    & 0.7 & 32.350 & 143.8 & 66.7 \\
    \textbf{4} 
    & 6.3 & 150.766 & 106.4 & 131.6
    & 5.0 & 187.040 & 113.4 & 143.1
    & 84 & {54.306} & 114.8 & 147.4 
    & 1.5 & \textbf{46.982} & 114.8 & 143.9 \\
    \textbf{5}
    & 6.3 & 119.159 & 147.1 & 71.3
    & 5.2 & 112.249 & 161.2 & 76.6
    & $>12$ h & $<52.052^\dagger$ & N/A & N/A
    & 2.4 & \textbf{48.265} & 162.2 & 76.6 \\
    \textbf{6}
    & 5.0 & 68.657 & 99.7 & 92.0 
    & 4.2 & 125.205 & 105.5 & 112.7
    & \multicolumn{4}{c|}{Timeout}
    & 1.1 & \textbf{29.482} & 111.3 & 104.6 \\
    \textbf{7}
    & 2.8 & 29.094 & 82.5 & 29.1
    & 3.1 & 34.213 & 83.3 & 32.6
    & \multicolumn{4}{c|}{Timeout}
    & 2.7 & \textbf{12.523} & 86.0 & 33.1 \\
    \textbf{8}
    & 3.0 & 20.065 & 79.3 & 22.0
    & 3.5 & 30.611 & 81.2 & 26.0
    & \multicolumn{4}{c|}{Timeout}
    & 3.2 & \textbf{9.001} & 82.4 & 23.0 \\
    \textbf{9}  
    & 7.0 & 174.850 & 105.5 & 130.3
    & 4.8 & 207.050 & 114.6 & 142.9
    & \multicolumn{4}{c|}{Timeout}
    & 4.1 & \textbf{43.401} & 124.3 & 145.2 \\
    \textbf{10}   
    & 4.9 & 121.066 & 104.1 & 84.5
    & 6.2 & 97.650 & 111.3 & 94.6
    & \multicolumn{4}{c|}{Timeout}
    & 7.5 & \textbf{26.912} & 125.5 & 91.5 \\
    \hline
    \textbf{Avg.} 
    & 4.813$\times$ & 3.012$\times$ & 0.922$\times$ & 0.926$\times$
    & 3.508$\times$ & 3.879$\times$ & 0.952$\times$ & 1.006$\times$
    & 14.204$\times$ & 1.048$\times$ & 0.995$\times$ & 1.002$\times$
    & \textbf{1}$\times$ & \textbf{1}$\times$ & \textbf{1}$\times$ & \textbf{1}$\times$\\ \hline
    \multicolumn{10}{c}{$^\dagger$ The binary of \texttt{SP-CP} only reports the temporary minimal total wirelength achieved at the time limit}
    \end{tabular}
    }
\end{table*}

\begin{table*}[htb]
    \centering
    \caption{Comparison under thermo-mechanical-aware optimization.
    Normalized metrics use \texttt{ATPlace2.5D} as baseline (1$\times$) for runtime and warpage, and \texttt{ATMPlace} for wirelength and temperature.}
    \label{tab:thermopt}
    \normalsize
    \resizebox{0.99\textwidth}{!}{
    \begin{tabular}{|c|cccc|cccc|cccc|cccc|}
    \hline
    \multirow{2}{*}{\textbf{Case}} 
    & \multicolumn{4}{c|}{\texttt{TAP-2.5D}~\cite{ma2021tap}} 
    & \multicolumn{4}{c|}{\texttt{TACPlace}~\cite{yu2025tacplace}}  
    & \multicolumn{4}{c|}{\texttt{ATPlace2.5D}~\cite{wang2024atplace2}} 
    & \multicolumn{4}{c|}{\texttt{ATMPlace (TM-aware)}} \\
    &RT/h&TWL/\SI{}{m}&$\text{T}_{\text{max}}/\SI{}{\degreeCelsius}$&$\text{Wpg}/\SI{}{\um}$
    &RT/h&TWL/\SI{}{m}&$\text{T}_{\text{max}}/\SI{}{\degreeCelsius}$&$\text{Wpg}/\SI{}{\um}$
    &RT/h&TWL/\SI{}{m}&$\text{T}_{\text{max}}/\SI{}{\degreeCelsius}$&$\text{Wpg}/\SI{}{\um}$
    &RT/h&TWL/\SI{}{m}&$\text{T}_{\text{max}}/\SI{}{\degreeCelsius}$&$\text{Wpg}/\SI{}{\um}$ \\ \hline
    \textbf{1}
    & 3.3 & 89.887 & 95.3 & 77.1
    & 5.3 & 62.040 & 100.6 & 85.6
    & 0.3 & 34.560 & 92.2 & 77.2
    & 0.4 & 25.900 & 94.6 & 76.5 \\
    \textbf{2}
    & 3.9 & 163.200 & 65.5 & 71.8
    & 5.0 & 97.773 & 66.5 & 95.2
    & 0.3 & 86.901 & 64.1 & 83.0
    & 0.3 & 56.807 & 64.9 & 84.3 \\
    \textbf{3}
    & 3.1 & 254.694 & 110.0 & 51.5
    & 4.5 & 114.597 & 139.7 & 69.3
    & 0.3 & 132.964 & 111.3 & 55.3
    & 0.4 & 219.062 & 107.3 & 34.1 \\
    \textbf{4}
    & 4.2 & 364.663 & 108.9 & 125.4
    & 4.8 & 180.313 & 111.4 & 147.0
    & 0.3 & 144.987 & 99.1 & 126.2
    & 0.4 & 150.416 & 100.6 & 128.3 \\
    \textbf{5}
    & 3.7 & 198.420 & 142.3 & 66.0
    & 6.5 & 116.413 & 163.8 & 73.7
    & 0.4 & 159.714 & 130.7 & 64.3
    & 0.5 & 117.776 & 123.2 & 60.0 \\
    \textbf{6}
    & 3.2 & 196.987 & 89.8 & 89.8
    & 6.8 & 121.524 & 102.2 & 105.3
    & 0.3 & 99.739 & 86.5 & 83.0
    & 0.5 & 82.236 & 83.3 & 79.3 \\
    \textbf{7}
    & 2.1 & 52.370 & 78.8 & 27.4
    & 6.8 & 33.602 & 82.7 & 32.6
    & 0.5 & 35.662 & 71.2 & 27.6
    & 0.5 & 36.119 & 70.7 & 27.6 \\
    \textbf{8}
    & 2.3 & 49.711 & 77.1 & 17.8
    & 7.4 & 32.767 & 79.5 & 24.2
    & 0.5 & 29.369 & 69.7 & 15.2
    & 0.6 & 23.680 & 69.8 & 16.4 \\
    \textbf{9}
    & 5.6 & 328.317 & 106.0 & 127.4
    & 8.4 & 194.374 & 118.5 & 144.0
    & 0.8 & 62.277 & 128.5 & 141.5
    & 0.9 & 70.717 & 118.6 & 140.0 \\
    \textbf{10}
    & 5.2 & 144.716 & 101.9 & 82.7
    & 8.4 & 106.826 & 120.1 & 97.1
    & 1.1 & 53.336 & 119.8 & 88.7
    & 1.2 & 59.547 & 125.2 & 90.8 \\
    \hline
    \textbf{Avg.}
    & 7.189 $\times$ & 2.312 $\times$ & 1.028 $\times$ & 1.048 $\times$ 
    & 12.203 $\times$ & 1.516 $\times$ & 1.132 $\times$ & 1.274 $\times$
    & 0.829 $\times$ & 1.101 $\times$ & 1.013 $\times$ & 1.063 $\times$
    & 1$\times$ & \textbf{1}$\times$ & \textbf{1}$\times$ & \textbf{1}$\times$ \\
    \hline
    \end{tabular}
    }
\end{table*}

\begin{figure*}[tbh]
\centering
    \begin{subfigure}[t]{0.24\linewidth}
        \centering
        \includegraphics[width=\linewidth]{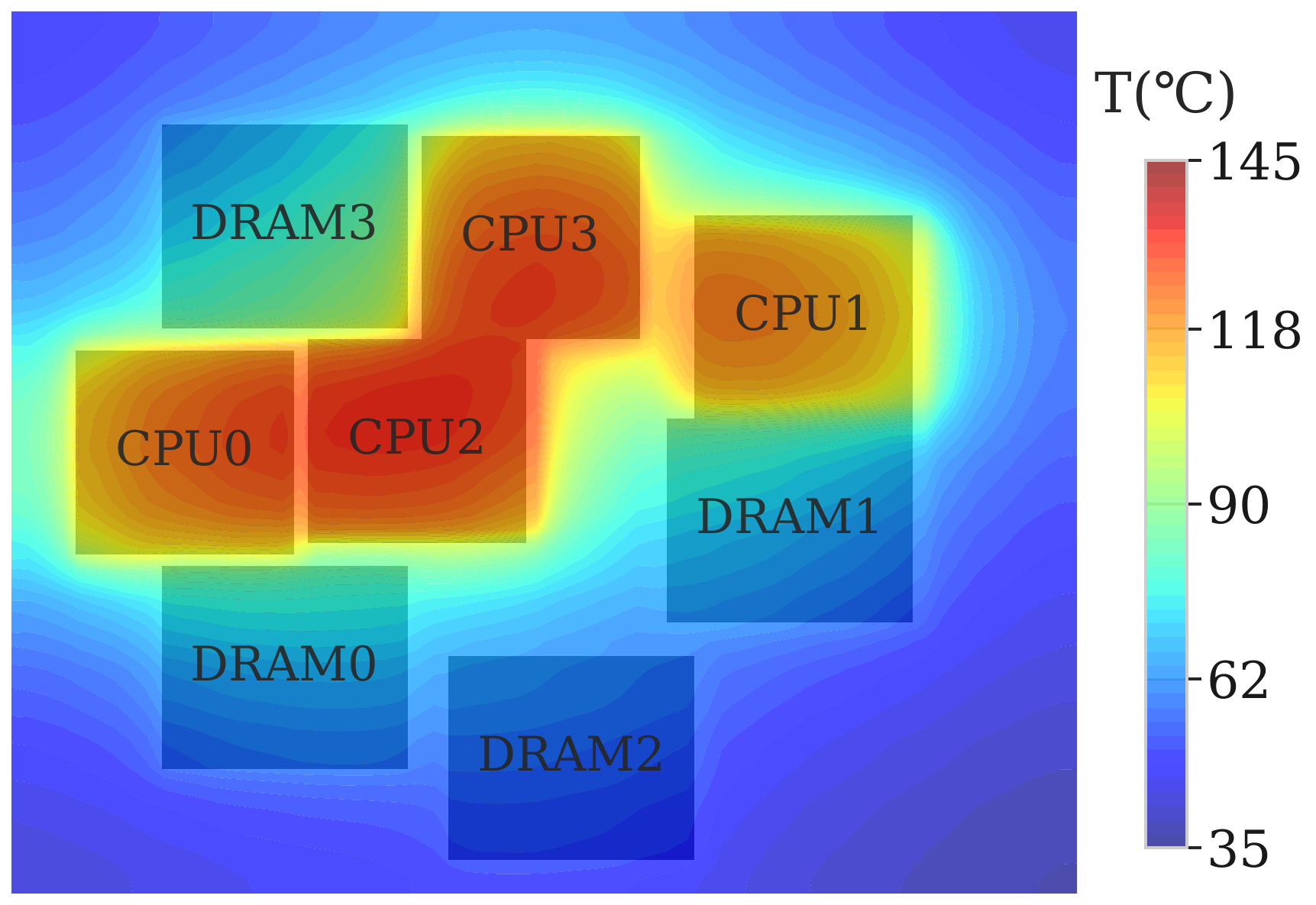}
        \caption{\texttt{TAP-2.5D} (WL-driven)}
    \end{subfigure}\hfill
    \begin{subfigure}[t]{0.24\linewidth}
        \centering
        \includegraphics[width=\linewidth]{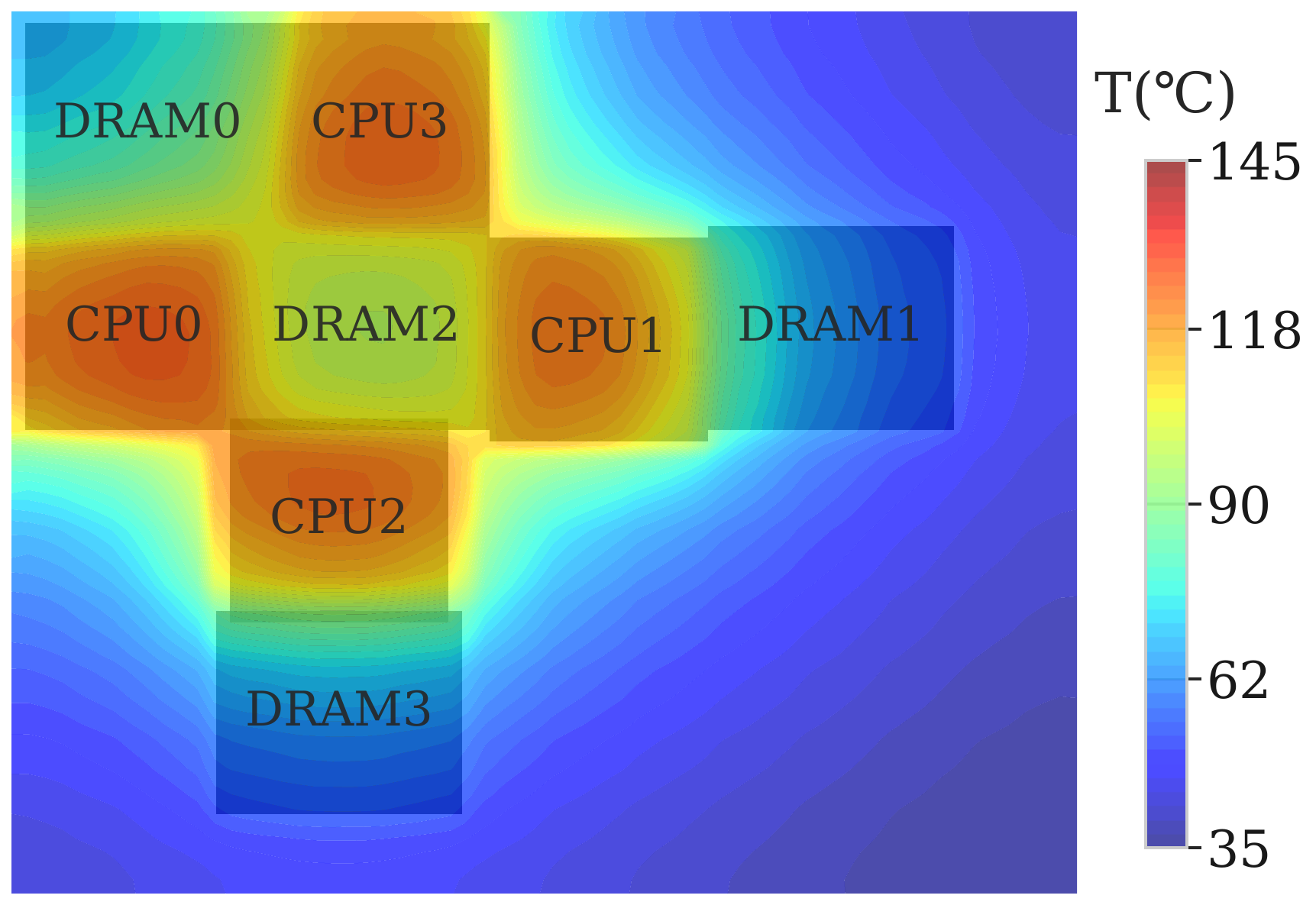}
        \caption{\texttt{TACPlace} (WL-driven)}
    \end{subfigure}\hfill
    \begin{subfigure}[t]{0.24\linewidth}
        \centering
        \includegraphics[width=\linewidth]{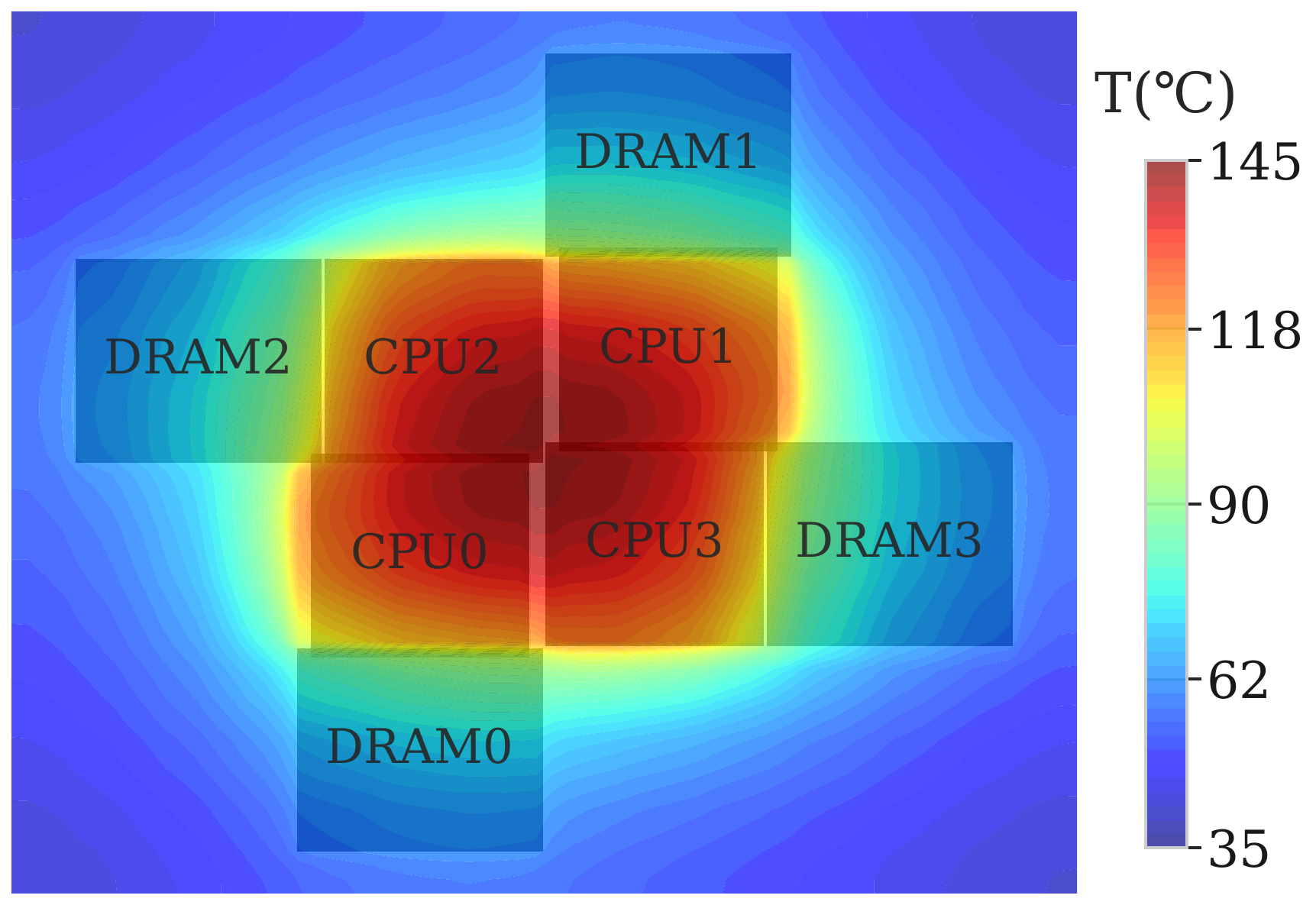}
        \caption{\texttt{SP-CP} (WL-driven)}
    \end{subfigure}\hfill
    \begin{subfigure}[t]{0.24\linewidth}
        \centering
        \includegraphics[width=\linewidth]{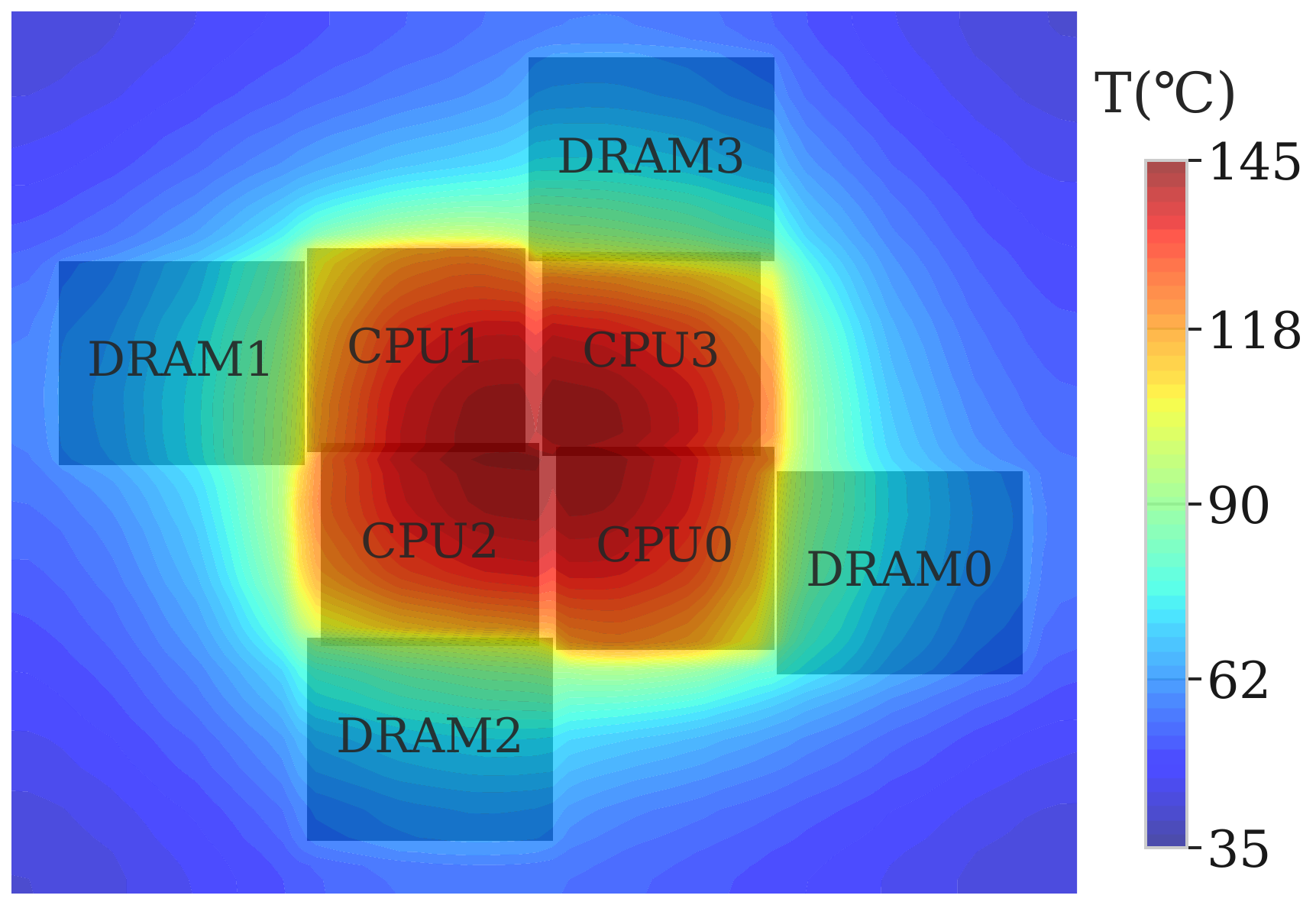}
        \caption{\texttt{ATMPlace} (WL-driven)}
    \end{subfigure}
\vspace{0.8em} 
    \begin{subfigure}[t]{0.24\linewidth}
        \centering
        \includegraphics[width=\linewidth]{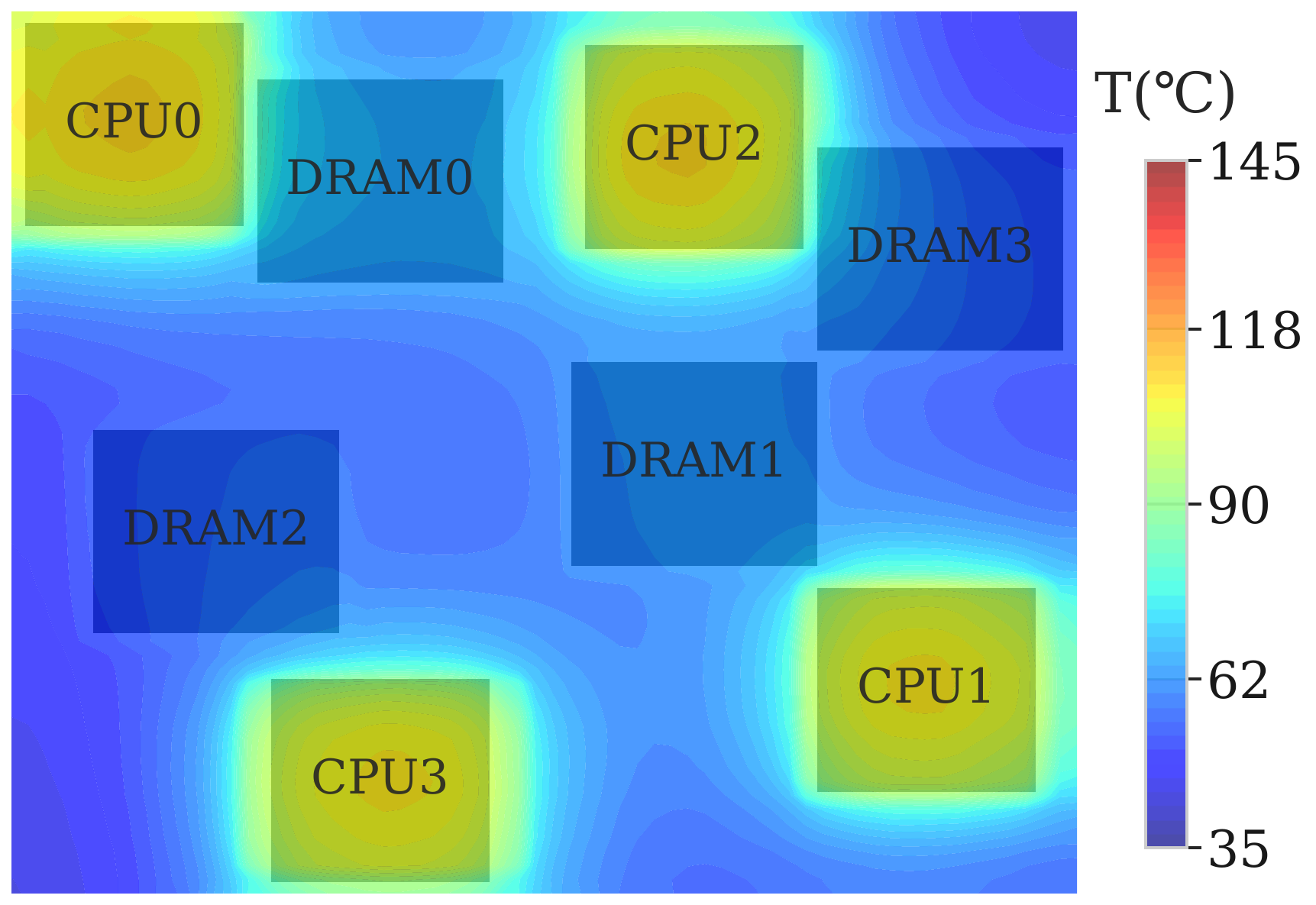}
        \caption{\texttt{TAP-2.5D} (TM-aware)}
    \end{subfigure}\hfill
    \begin{subfigure}[t]{0.24\linewidth}
        \centering
        \includegraphics[width=\linewidth]{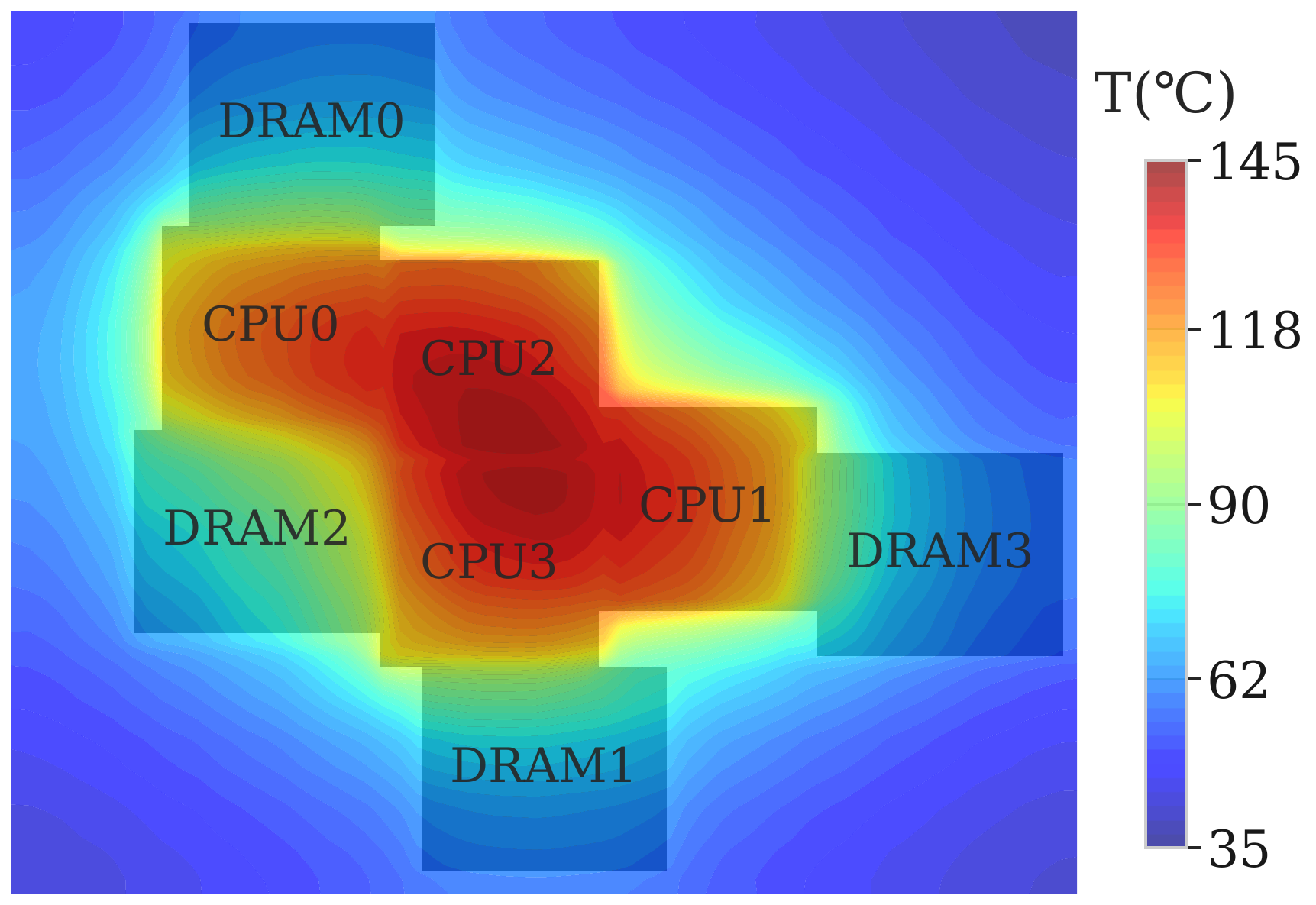}
        \caption{\texttt{TACPlace} (TM-aware)}
    \end{subfigure}\hfill
    \begin{subfigure}[t]{0.24\linewidth}
        \centering
        \includegraphics[width=\linewidth]{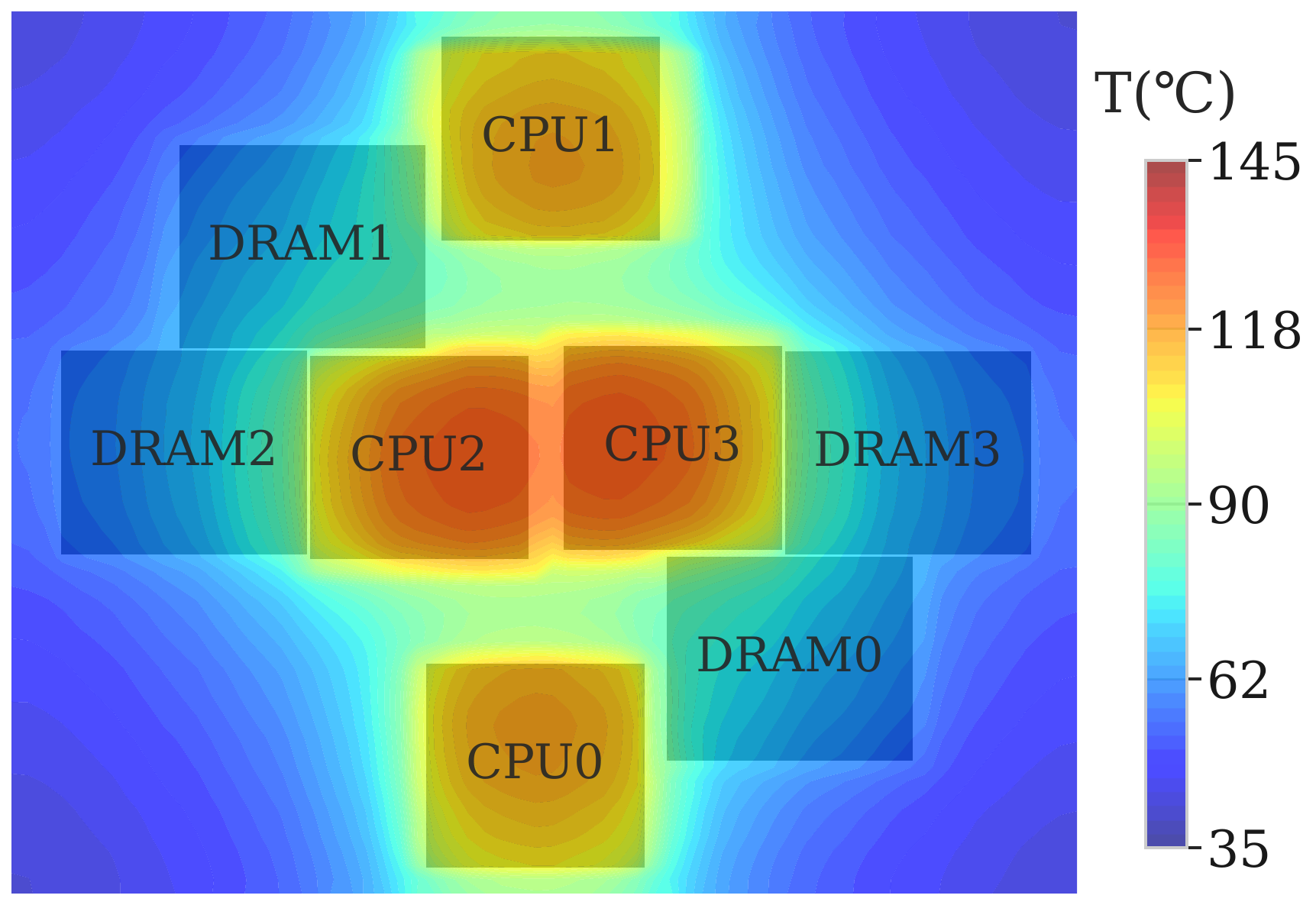}
        \caption{\texttt{ATPlace2.5D} (T-only)}
    \end{subfigure}\hfill
    \begin{subfigure}[t]{0.24\linewidth}
        \centering
        \includegraphics[width=\linewidth]{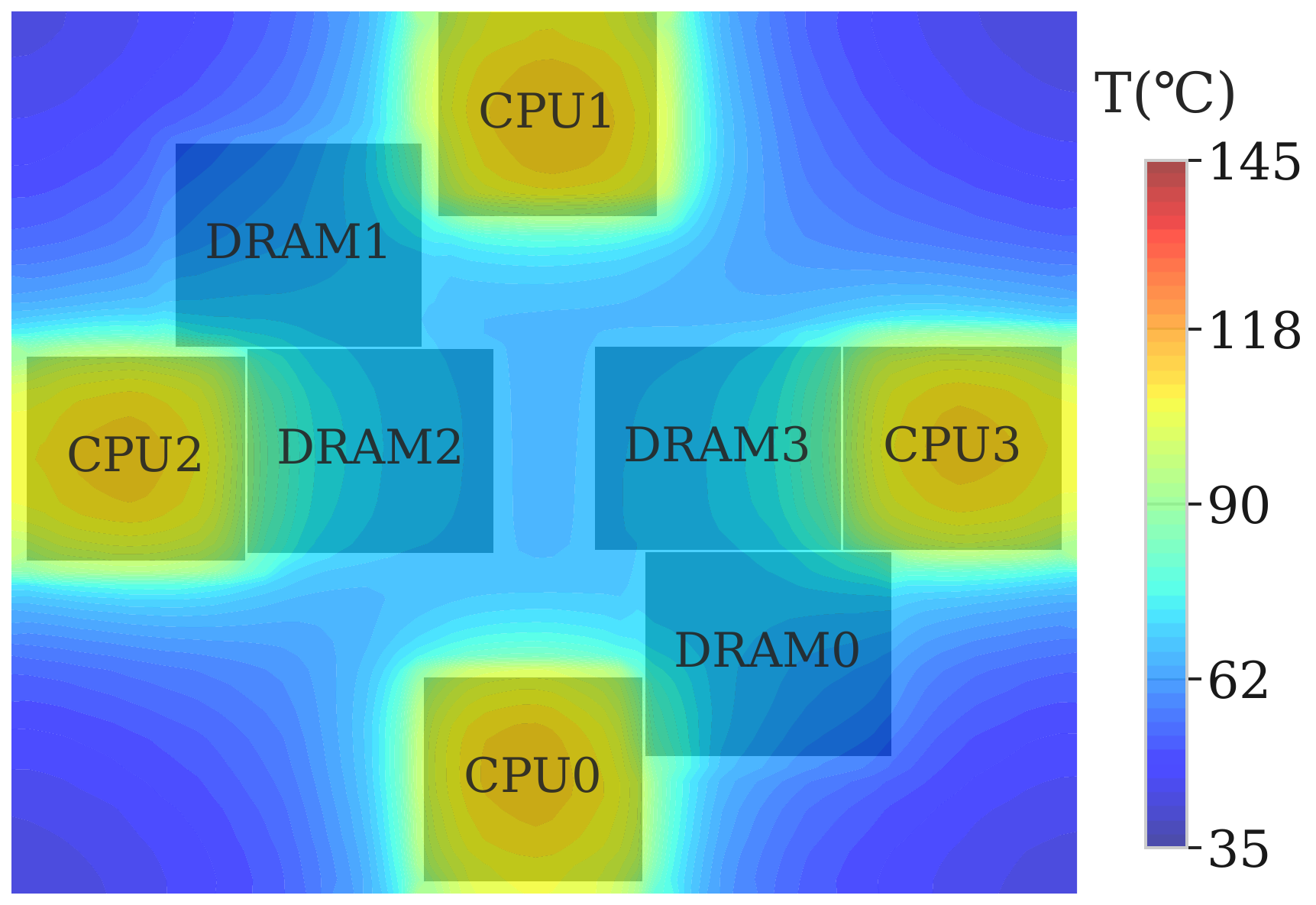}
        \caption{\texttt{ATMPlace} (TM-aware)}
    \end{subfigure}
\caption{
    Placement layouts for Case~3. Upper row: wirelength-driven optimization. Lower row: thermo-mechanical-aware optimization. 
    CPU and DRAM dimensions: $9 \times 8~\mathrm{mm}^2$ and $9 \times 9~\mathrm{mm}^2$, respectively.
}
\label{fig:layouts}
\end{figure*}

\subsection{Placement Results}
\label{sec:placement}

We evaluate four representative placement methodologies: 
\begin{itemize}
    \item \texttt{TAP-2.5D}~\cite{ma2021tap}: Simulated Annealing (SA)-based, open-sourced and adapted to our benchmark;
    \item \texttt{TACPlace}~\cite{yu2025tacplace}: GA/PSO hybrid, representing evolutionary approaches, realized and adapted for thermo-mechanical applications;
    \item \texttt{SP-CP}~\cite{chiou2023chiplet}: enumeration-based (times out for cases with more than 12 chiplets), wirelength-driven binary tested only;
    \item \texttt{ATPlace2.5D}~\cite{wang2024atplace2}: our prior work, thermal-aware only;
    \item \texttt{ATMPlace}: the proposed framework, jointly optimizing thermal and mechanical warpage.
\end{itemize}
Reinforcement learning (RL)-based methods~\cite{duan2023rlplanner,deng2024chiplet} are excluded due to unavailability of executables; reported results suggest their performance is comparable to or slightly better than \texttt{TAP-2.5D}. The method of \cite{osmolovskyi2018optimal} fails on our benchmarks due to scalability limitations.

\subsubsection{Wirelength-Driven Optimization}
\label{sec:wl_opt_analysis}

For completeness and to isolate the impact of physical-aware objectives, we first evaluate pure wirelength-driven (WL-driven) optimization. 
As shown in Table~\ref{tab:wlopt}, \texttt{ATMPlace} achieves state-of-the-art wirelength minimization: it delivers the shortest TWL in 8 out of 10 cases, and ties \texttt{SP-CP} in Case~1 and 3. 
The average TWL is $1.00\times$ (baseline), compared to $3.0\times$ and $3.9\times$ for \texttt{TAP-2.5D} and \texttt{TACPlace}, respectively. Even in Case~5—the largest enumerated instance (\texttt{SP-CP} times out with $< \SI{52.052}{m}$ TWL)—\texttt{ATMPlace} obtains TWL $= \SI{48.265}{m}$, which is $>2.4\times$ shorter than \texttt{TAP-2.5D} and $>2.3\times$ shorter than \texttt{TACPlace}.
Moreover, the runtime advantage of \texttt{ATMPlace} persists in WL mode: with average runtime $= \SI{2.3}{min}$, it is $14.2\times$ faster than \texttt{SP-CP} (on solvable cases) and $>3\times$ faster than heuristic methods. The enumeration approach, while yielding marginally better TWL (e.g., in Case~1), fails to scale beyond 12 chiplets and produces thermally hazardous layouts. 

This superiority arises from the synergy of three components: (i) the MILP-based initialization, which yields a high-quality seed layout with orientation-aware bump alignment; (ii) the smooth orientation parametrization, which enables gradient-based refinement of discrete rotations without combinatorial explosion; and (iii) the adaptive density-weight scheduling ~\eqref{eq:overflow}, which dynamically balances HPWL reduction against placement overflow—preventing premature congestion that traps heuristic methods.

However, WL-driven optimization incurs significant physical penalties. For instance, in Case~3, \texttt{ATMPlace} achieves TWL $= \SI{32.350}{m}$ (vs. $\SI{134.073}{m}$ in TM mode) but suffers $T_{\max} = \SI{143.8}{\degreeCelsius}$—an increase of \SI{36.5}{\degreeCelsius} over the TM-aware variant—and warpage rises from $\SI{34.1}{\um}$ to $\SI{66.7}{\um}$. 
Similar trends hold universally: minimizing TWL alone encourages tight clustering of connected chiplets (e.g., GPUs with adjacent HBMs), intensifying local hotspots and bending moments. 
This confirms the \textit{necessity} of multi-physics co-optimization: wirelength-driven placement, while optimal in isolation, is physically infeasible for advanced 2.5D-ICs under stringent thermo-mechanical reliability constraints.

\begin{figure*}[htb]
    \centering
    \begin{subfigure}[t]{0.48\linewidth}
        \centering
        \includegraphics[width=\textwidth]{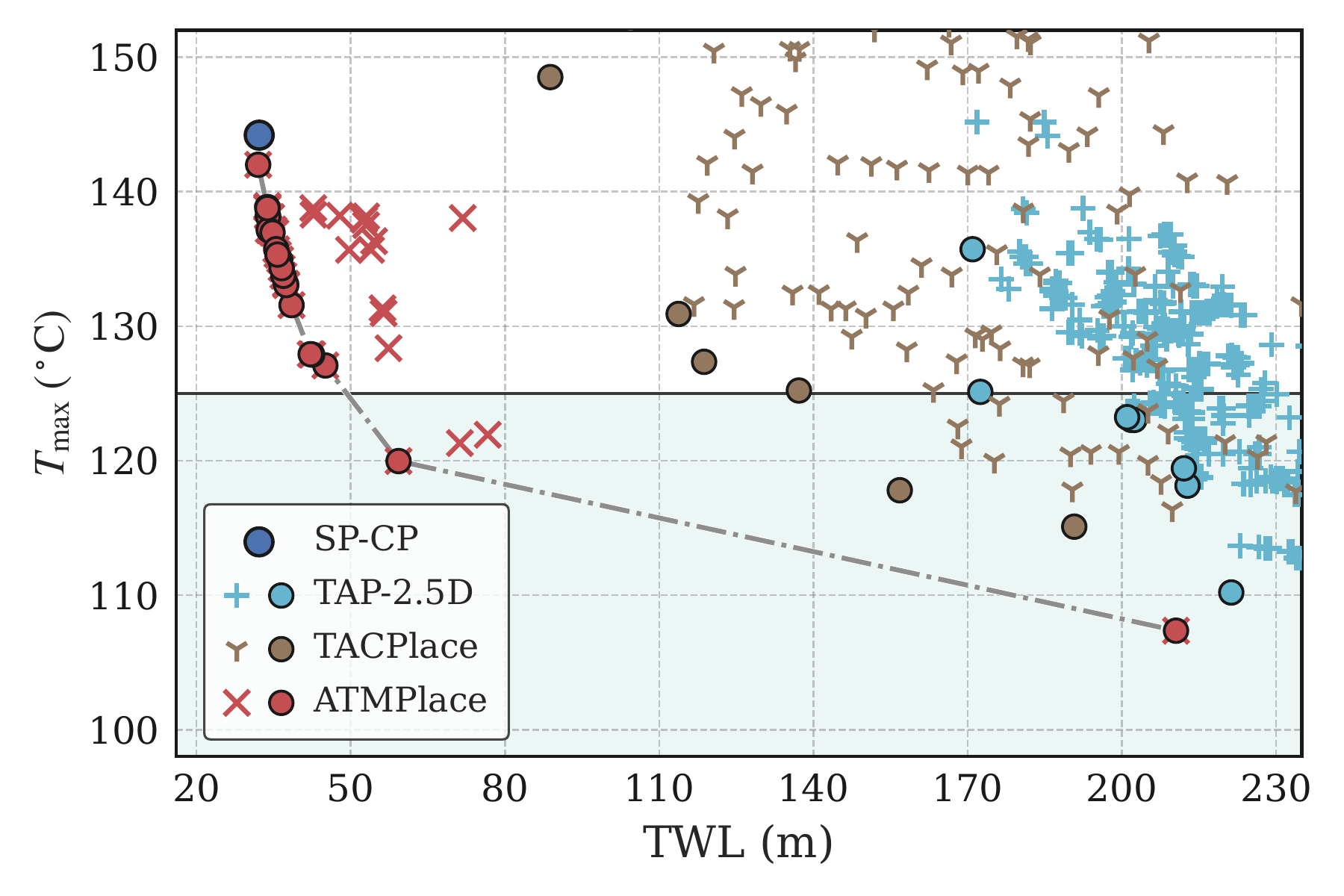}\\
        \caption{Case~3}
    \end{subfigure}
    \hfill
    \begin{subfigure}[t]{0.48\linewidth}
        \centering
        \includegraphics[width=\textwidth]{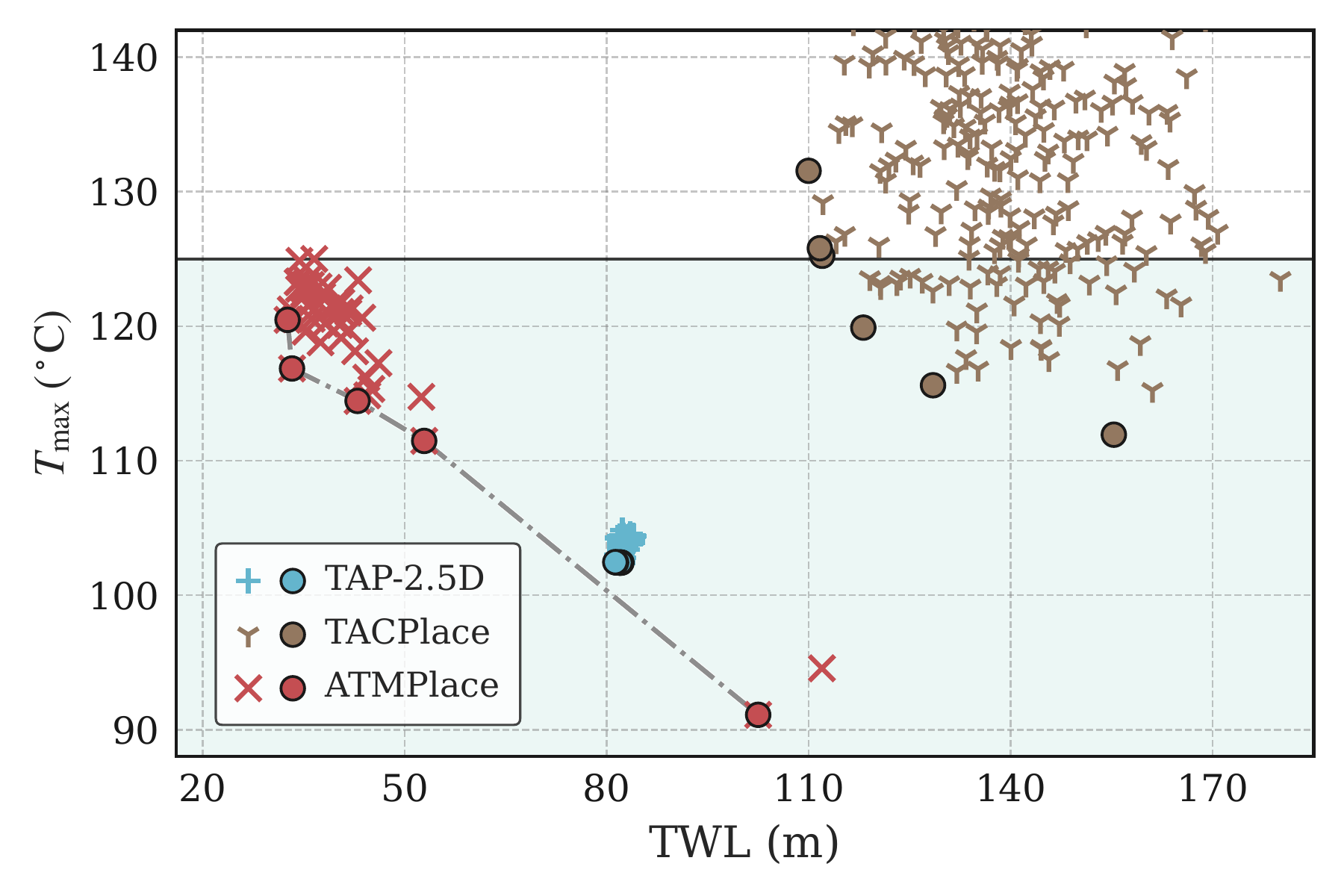}\\
        \caption{Case~10}
    \end{subfigure}
    \caption{Multi-objective trade-offs between maximum temperature ($T_{\max}$) and total wirelength (TWL) for Case~3 (left) and Case~10 (right). Solid circles represent the Pareto-optimal solutions obtained by each method. The dashed line indicates an approximate Pareto frontier constructed from the best-known solutions across all methods.}
    \label{fig:pareto}
\end{figure*}

\subsubsection{Thermo-Mechanical-Aware Optimization}
\label{sec:tm_opt}
\begin{figure}[tbh]
    \centering
    \begin{subfigure}[t]{0.49\linewidth}
        \centering
        \includegraphics[width=1\linewidth]{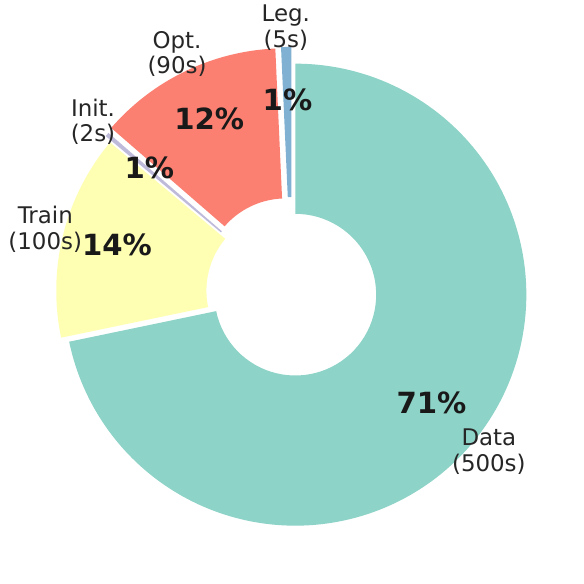}\\
        \caption{Case 3}
    \end{subfigure}
    \hfill
    \begin{subfigure}[t]{0.48\linewidth}
        \centering
        \includegraphics[width=1\linewidth]{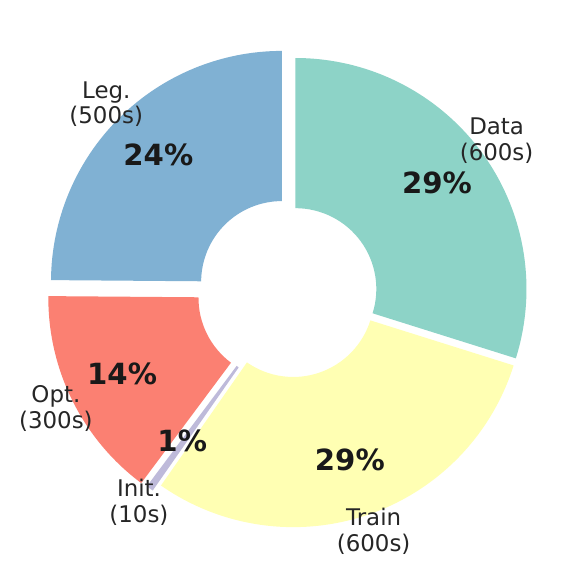}\\
        \caption{Case 10}
    \end{subfigure}
    \caption{Runtime breakdown of \texttt{ATMPlace} in thermo-mechanical-aware optimization. `Data.': dataset generation; `Train.': compact model training; `Init.': initialization; `Opt.': optimization; `Leg.': legalization. Percentages indicate the share of the total runtime.}
    \label{fig:runtime}
\end{figure}

We summarize the thermo-mechanical-aware  placement results in Table~\ref{tab:thermopt}. Across all ten benchmarks, \texttt{ATMPlace} consistently achieves the best trade-off among wirelength, thermal, and mechanical objectives on average.
The average gains are decisive: normalized to \texttt{ATMPlace}, \texttt{TAP-2.5D} and \texttt{TACPlace} incur 2.5$\times$ and 1.5$\times$ longer TWL, 3\%–13\% higher $T_{\max}$, and 5\%–27\% larger warpage—at 7$\times$–12$\times$ runtime overhead. 
Runtime averages \SI{31}{min}—7.5$\times$ faster than \texttt{TAP-2.5D}—with per-run optimization costing only 20–40\% of total time after one-time model training (Fig.~\ref{fig:runtime}). 

Besides, \texttt{ATMPlace} yields the lowest $T_{\max}$ in 6 cases, the shortest TWL in 5 cases, and the smallest warpage in 5 cases. 
Notably, it simultaneously minimizes all three metrics in Case~5 and 6. In Case~5, \texttt{ATMPlace} reduces $T_{\max}$ by \SI{19.1}{\degreeCelsius} and warpage by \SI{6.3}{\um} over \texttt{TAP-2.5D}, while shortening TWL by 40\%. 
Even against our thermal-optimization-only predecessor \texttt{ATPlace2.5D}, \texttt{ATMPlace} lowers warpage by up to 60\% (Case~3) without degrading $T_{\max}$. 
Thus, \texttt{ATMPlace} is the only method achieving scalable, high-quality TM co-optimization.
This stems from three key features:  
(i) \textit{Joint differentiable modeling} of thermal and warpage fields, enabling gradient-guided symmetry enforcement (e.g., mirrored placement of high-power dies);  
(ii) \textit{Warpage-aware regularization} via ~\eqref{eq:warpage_compact}, which suppresses torque-inducing layouts missed by thermal-only optimization;  
(iii) \textit{Decoupled step-size control} for position/orientation (Algorithm~\ref{alg:cgd}), preventing angular updates from destabilizing mechanical equilibrium.

\subsubsection{Placement Layouts}

Figure~\ref{fig:layouts} illustrates the placement outcomes for Case~3 under wirelength-driven (WL-driven) and thermo-mechanical-aware (TM-aware) optimizations, revealing distinct layout strategies dictated by objective formulation. 

Under wirelength-driven optimization, the enumeration-based \texttt{SP-CP} (Fig.~\ref{fig:layouts}(c)) achieves the global WL minimum by aggregating the four CPUs at the center and arranging DRAMs radially—a layout that minimizes interconnect distances but concentrates thermal power density. 
Remarkably, \texttt{ATMPlace} (Fig.~\ref{fig:layouts}(d)) reproduces this near-optimal topology with only $0.5\%$ WL overhead, demonstrating the efficacy of its gradient-based refinement from a high-quality initial solution. 
In contrast, heuristic methods exhibit structural inefficiencies: \texttt{TAP-2.5D} (Fig.~\ref{fig:layouts}(a)) splits the CPUs, increasing inter-CPU and CPU–DRAM wirelength, while \texttt{TACPlace} (Fig.~\ref{fig:layouts}(b)) yields an elongated, loosely packed configuration—both indicative of premature convergence in high-dimensional discrete spaces.

In thermo-mechanical-aware optimization, physical constraints dominate layout morphology. \texttt{TAP-2.5D} (Fig.~\ref{fig:layouts}(e)) adopts conservative clustering, resulting in low warpage ($\SI{52}{\um}$) but poor thermal performance ($T_{\max} = \SI{110}{\degreeCelsius}$) and excessive wirelength (+124\% over \texttt{ATMPlace}). \texttt{TACPlace} (Fig.~\ref{fig:layouts}(f)) partially disperses CPUs yet retains asymmetry, leading to  the highest temperature and warpage ($\SI{140}{\degreeCelsius}, \SI{69}{\um}$) among TM-aware methods. 
The thermal-only \texttt{ATPlace2.5D} (Fig.~\ref{fig:layouts}(g)) improves temperature ($T_{\max} = \SI{111}{\degreeCelsius}$) via CPU separation but overlooks mechanical coupling—its staggered DRAM placement induces torque-like bending, reflected in elevated warpage ($\SI{55}{\um}$). 
Only \texttt{ATMPlace} (Fig.~\ref{fig:layouts}(h)) achieves full co-optimization: CPUs are symmetrically positioned, DRAMs are mirrored across both axes, and inter-chiplet spacing is uniformly enlarged in high-gradient zones. 
This layout delivers balanced $T_{\max} = \SI{107}{\degreeCelsius}$ and $Wpg = \SI{34}{\um}$—matching the thermal performance of \texttt{ATPlace2.5D} while ensuring mechanical robustness through spatial symmetry enforced by the warpage-aware gradient term.

\begin{figure}[hbt]
    \centering
    \includegraphics[width=0.95\linewidth]{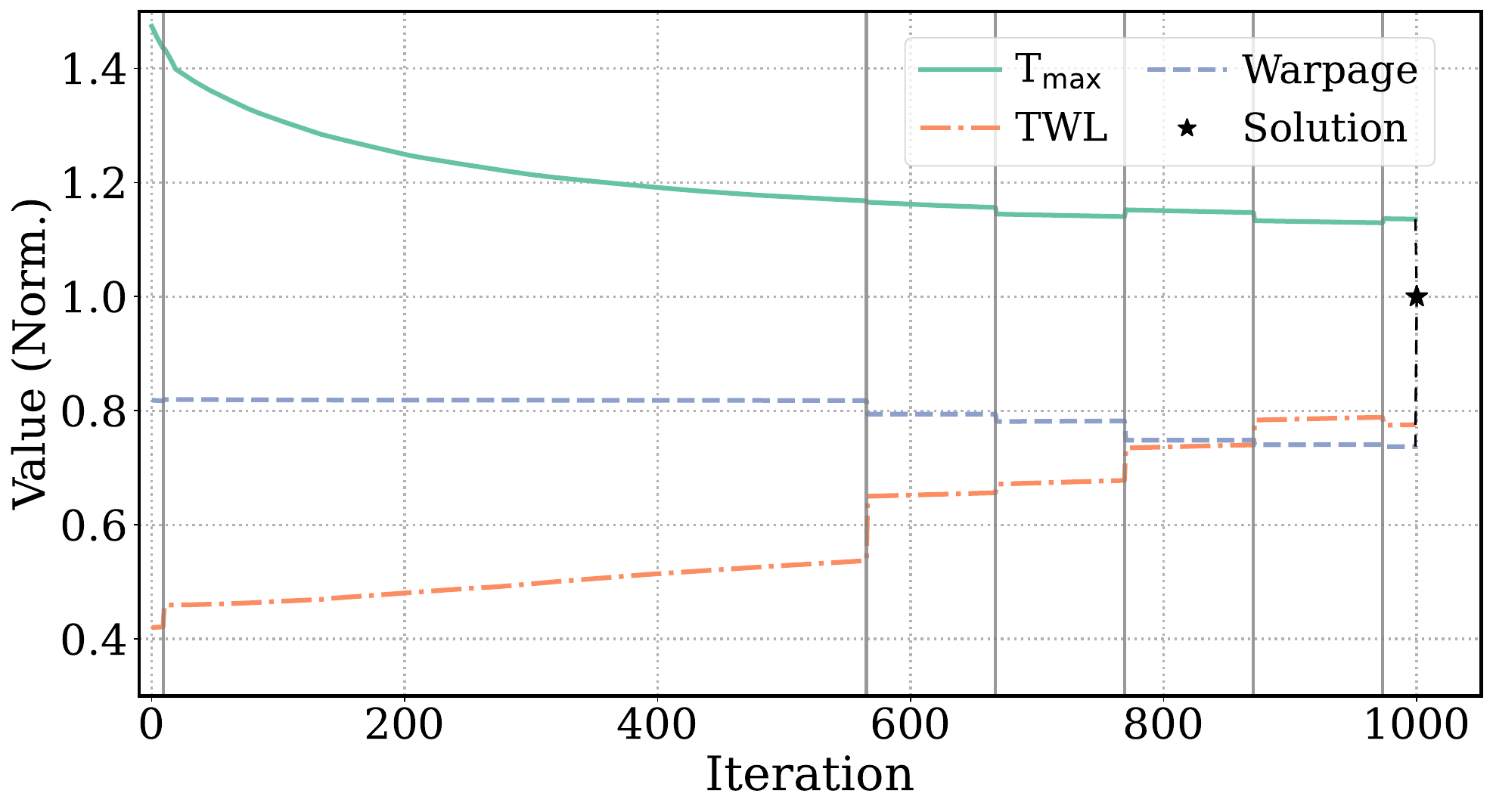}
    \caption{Normalized convergence trajectories of $T_{\max}$, wirelength (TWL), and warpage under thermo-mechanical-aware optimization (Case~3). Curves show evolution before the final solution (marked by $\star$); the vertical lines denote entropy-driven noise injections. All metrics are normalized to their final values.}
    \label{fig:converge}
\end{figure}

\subsubsection{Multi-Objective Optimization}

Figure~\ref{fig:pareto} quantifies the trade-offs between $T_{\max}$ and TWL for Case~3 (left) and Case~10 (right). 
\texttt{ATMPlace} generates a dense, continuous Pareto front spanning from wirelength-optimal to thermo-mechanical-optimal, demonstrating full control over the design space.

In Case~3, \texttt{ATMPlace} dominates all baselines: its solutions form a tight, low-$T_{\max}$ cluster below $\SI{125}{\degreeCelsius}$, while \texttt{TAP-2.5D} and \texttt{TACPlace} are confined to higher wirelength region. The enumeration-based \texttt{SP-CP} achieves minimal TWL but suffers extreme thermal penalty ($T_{\max} > \SI{140}{\degreeCelsius}$).
For the large-scale Case~10, \texttt{ATMPlace} again establishes the frontier, while \texttt{TAP-2.5D} clusters, and  \texttt{TACPlace} fail to enter the optimal region. This highlights \texttt{ATMPlace}'s scalability and ability to navigate complex multi-objective landscapes where heuristic methods stagnate.
To conclude, \texttt{ATMPlace} establishes the Pareto frontier: it dominates all baselines in the three-dimensional objective space of (TWL, temperature, warpage).


\subsubsection{Runtime Analysis}
\label{sec:runtime}
Figure~\ref{fig:runtime} details the runtime breakdown for Case~3 and Case~10. Model construction (dataset generation + training) dominates (\SI{85}{\%} for Case~3, \SI{58}{\%} for Case~10), but is a one-time cost amortized over multiple runs (e.g., Pareto exploration). Initialization is negligible ($<\SI{10}{\second}$). The core optimization and legalization stages consume \SI{10}{\%}–\SI{40}{\%} of total time, reflecting their computational complexity.

Figure~\ref{fig:converge} shows the normalized convergence trajectories for Case~3 under thermo-mechanical-aware optimization. 
As high-power chiplets disperse to mitigate thermal hotspots, $T_{\max}$ drops sharply in early iterations (from 1.4$\times$ to $\sim$1.15$\times$ its final value), while TWL initially rises—reflecting the trade-off between thermal relief and interconnect cost—before stabilizing near its final value. 
Warpage exhibits a gradual decline, indicating progressive mechanical stress reduction throughout the optimization.
Periodic entropy-driven noise injections (dashed vertical lines) perturb the layout when the optimization gets stuck, enabling escape from local minima and yielding incremental improvements across all three objectives. Notably, after each injection, the system converges faster to a lower-cost state, demonstrating that the perturbation enhances exploration without destabilizing the solution trajectory.

The robustness of our approach is further evidenced by the smooth, deterministic convergence profile: unlike stochastic methods (e.g., SA, GA), where solution quality varies significantly with random seeds—posing a reliability risk in the design flow—our method consistently reaches the final solution ($\star$). The final solution (at iteration 1000) represents a balanced Pareto-optimal trade-off among thermal, mechanical, and wiring constraints.

\section{Conclusion}\label{sec:Conclusion}
As the semiconductor industry increasingly embraces 2.5D-ICs to overcome the limitations of traditional scaling, the design of large-scale, heterogeneous chiplet-based systems has become both urgent and inevitable.  
This paper presents \texttt{ATMPlace}, a scalable analytical placement framework for large-scale 2.5D-ICs that jointly optimizes wirelength, thermal integrity, and mechanical reliability—addressing a critical gap in existing methodologies.
Our physics-based compact thermal-mechanical model attains a mean correlation above {0.97} and accelerates simulation by {over $8000\times$} compared to numerical solvers, validating its fidelity for gradient-based optimization.  
By integrating orientation-aware analytical placement with fast compact models, \texttt{ATMPlace} achieves high-quality layouts at unprecedented speed: it delivers up to {146\%} and {52\%} better total wirelength, {3\%} and {13\%} lower peak temperature, and {5\%} and {27\%} better warpage than \texttt{TAP-2.5D} and \texttt{TACPlace}, respectively, under thermo-mechanical-aware constraints.  
In wirelength-driven mode, it matches the near-optimal \texttt{SP-CP} within 1\% deviation.

Looking ahead, we plan to extend the framework to support co-design with interposer routing by incorporating congestion and signal integrity into the objective, enabling joint placement-routing optimization.  
We believe \texttt{ATMPlace} will lay a critical foundation for scalable, reliable, and automated design of next-generation 2.5D-IC systems.

\section{Acknowledgment}
The authors thank Xueqing Li, Yifan Chen, Xinming Wei for insightful discussions. This paper used GPT-5 to generate the thermal deformation schematic figure in Sec.~I.

{
\small
\bibliographystyle{IEEEtranN}
\bibliography{./Ref/ATPlace}
}

\end{document}